\documentclass[hidelinks,onefignum,onetabnum]{siamart250211}

\usepackage{lipsum}
\usepackage{amsfonts}
\usepackage{comment}
\usepackage{graphicx}
\usepackage{epstopdf}
\usepackage{algorithm}
\usepackage{algpseudocode}
\usepackage{multirow}
\usepackage{array}
\usepackage{float}
\usepackage{tabularx}
\usepackage{booktabs}
\ifpdf
  \DeclareGraphicsExtensions{.eps,.pdf,.png,.jpg}
\else
  \DeclareGraphicsExtensions{.eps}
\fi

\usepackage{amssymb}

\usepackage{subcaption}

 \usepackage{amsopn}

\DeclareMathOperator{\ReLU}{ReLU}

\newcommand{\bcF}{\boldsymbol{\mathcal{F}}}
\newcommand{\bka}{\boldsymbol{\kappa}}
\newcommand{\brh}{\boldsymbol{\rho}}
\newcommand{\bth}{\boldsymbol{\theta}}
\newcommand{\bu}{\boldsymbol{u}}
\newcommand{\bx}{\boldsymbol{x}}
\newcommand{\bz}{\boldsymbol{z}}

\newsiamremark{remark}{Remark}
\newsiamremark{hypothesis}{Hypothesis}
\crefname{hypothesis}{Hypothesis}{Hypotheses}
\newsiamthm{claim}{Claim}
\newsiamremark{fact}{Fact}
\crefname{fact}{Fact}{Facts}

\headers{FM-tfPINN for tTFCP Systems}{Kumar et al.}

\title{A Fractional-Memory Physics-Informed Neural Network with Fast History Compression for Tempered Fractional Coupled Phase-Field Systems}

\author{
Shubham Kumar\thanks{Department of Mathematical Sciences, IIT (BHU), Varanasi 221005, India
  (\email{shubhamkumar.rs.mat24@itbhu.ac.in},
   \email{himanshukrdwivedi.rs.mat21@itbhu.ac.in},
   \email{rajeev.apm@iitbhu.ac.in}).}
\and
Himanshu Kumar Dwivedi\footnotemark[1]
\and
Matthias Ehrhardt\thanks{Corresponding author. IMACM, School of Mathematics and Natural Sciences,
University of Wuppertal, Germany
  (\email{ehrhardt@uni-wuppertal.de}).}
\and
Rajeev\footnotemark[1]
}

\ifpdf
\hypersetup{
  pdftitle={FM-tfPINN: A Fractional-Memory Generated Physics-Informed Neural Network for Tempered Fractional Coupled Phase-Field Systems},
  pdfauthor={S. Kumar, H.K. Dwivedi, and Rajeev}
}
\fi

\begin{document}

\maketitle

\begin{abstract}
\textit{Tempered time-fractional coupled phase-field} (tTFCP) systems are used to model interfacial phenomena involving memory-dependent transport and relaxation mechanisms. 
Numerical solutions to these systems are challenging due to the simultaneous presence of nonlocal temporal operators, weak initial singularities, moving diffuse interfaces, and strongly coupled multiphysics dynamics. 
In this work, we introduce FM-tfPINN (fractional-memory physics-informed neural network), which is used for forward simulation and inverse parameter identification in tempered fractional coupled phase-field systems. 
Unlike conventional fractional PINNs, which enforce memory effects solely through residual constraints, our framework incorporates tempered fractional memory directly into the neural representation via latent memory-source functions and a tempered fractional integral operator. 

We develop a fast shifted residual formulation based on graded temporal meshes and \textit{sum-of-exponentials} (SOE) history compression to efficiently evaluate the tempered fractional operators. 
This framework combines interface-aware and residual-adaptive collocation strategies, improving resolution near evolving diffuse interfaces. 
A unified, physics-informed loss formulation allows for the forward prediction and inverse recovery of unknown physical parameters from sparse observations.

We assess the proposed method on a class of tempered fractional corrosion phase-field models, including one-dimensional corrosion-front propagation, activation- and diffusion-controlled regimes, two-dimensional pitting corrosion, and inverse mobility identification problems. 
The numerical results demonstrate the accurate recovery of coupled phase and concentration fields, the robust prediction of physically relevant interface diagnostics, and the reliable estimation of parameters from limited data. 
The FM-tfPINN framework efficiently and consistently learns and solves nonlinear tempered fractional phase-field systems, offering a general paradigm for physics-informed computation of nonlocal evolution equations.

\end{abstract}



\section{Introduction}
Phase-field modeling is a powerful approach for describing interfacial evolution in complex multiphysics systems. 
It replaces sharp moving boundaries with smoothly varying order parameters, thus eliminating avoiding the need for explicit front tracking
\cite{BoettingerWarrenBeckermannKarma2002, Caginalp1986, Chen2002PhaseField, ProvatasElder2010}.
This diffuse-interface viewpoint is widely used in reaction-diffusion dynamics, phase separation, microstructure evolution, and related interfacial processes \cite{LiuChengWangZhao2018, TangYuZhou2019}.
In corrosion science, phase-field formulations are particularly attractive because the metal surface degradation is governed not only by the motion of the metal-electrolyte interface, the redistribution of chemical concentrations, and the interaction between interfacial kinetics and transport mechanisms \cite{AnsariXiaoHuLiLuoShi2018, ChenLucariniMaChenCui2025PFPINNs, MaiSoghratiBuchheit2016}. 
Consequently, coupled phase-field models naturally lend themselves to studying corrosion-front propagation, pitting morphology, and concentration-driven interfacial evolution in a unified framework.
Classical integer-order phase-field models may be insufficient when the underlying dynamics exhibit history-dependent transport, anomalous relaxation, or memory-influenced interfacial kinetics. 
Time-fractional models naturally describe such effects and have been widely used to represent nonlocal temporal behavior in diffusion and transport processes \cite{KilbasSrivastavaTrujillo2006, MetzlerKlafter2000, Podlubny1999}. 
In many practical settings, however, purely power-law memory may overemphasize very long histories. 
Tempered fractional models provide a more flexible description by retaining fractional memory effects over relevant time scales while introducing memory attenuation at longer times \cite{BaeumerMeerschaert2010, MeerschaertZhangBaeumer2008, SabzikarMeerschaertChen2015}.
This makes time-fractional modeling well-suited for corrosion-driven phase-field dynamics, in which interfacial motion and concentration redistribution depend on accumulated history and gradual memory decay.

The numerical treatment of tTFCP systems is challenging due to 
several sources of complexity.
The temporal derivative is nonlocal, requiring an accurate representation of the accumulated history. 
Meanwhile, the solution may exhibit limited regularity near the initial time \cite{JinLazarovZhou2016, McLean2010,dwivedi_ehrhardt, StynesORiordanGracia2017}. 
Next, the phase-field and concentration variables are strongly coupled, and the dominant dynamics are often localized near thin diffuse-interface regions.
Standard time-marching solvers thus 
require careful temporal grading, efficient history compression, and sufficient spatial resolution near moving interfaces. 
Although fast convolution techniques based on sum-of-exponentials (SOE) approximations have significantly reduced the cost of fractional-memory evaluation \cite{Arnold2003, DwivediRajeevZeng2026JSC, JiangZhangZhangZhang2017}, the construction of flexible, learning-based solvers for coupled, tempered, fractional, phase-field dynamics remains challenging.

\textit{Physics-informed neural networks} (PINNs) are a versatile, mesh-free learning paradigm for solving forward and inverse problems governed by differential equations. They do this by incorporating underlying physical laws directly into the optimization loss \cite{RaissiPerdikarisKarniadakis2019PINNs, KarniadakisKevrekidisLuPerdikarisWangYang2021}.
Efficient realization of residual-driven neural discretizations is facilitated by automatic differentiation, which enables systematic evaluation of high-order differential operators within contemporary machine learning platforms \cite{BaydinPearlmutterRadulSiskind2017, PaszkeGrossChintalaEtAl2017}.
Several neural approaches have also been proposed for PDEs exhibiting oscillatory behavior or requiring the preservation of underlying structural properties \cite{CaiLiLiu2020PhaseDNN}. 
Extending PINNs to fractional models enables data-driven approximation and parameter estimation for several classes of nonlocal problems \cite{PangLuKarniadakis2019,RenMengLiuHouYu2023,ShiLiuYang2025}. 
However, direct PINN and fractional PINN (fPINN) formulations may encounter challenges with tTFCP systems because the memory structure is typically imposed through the residual loss,
whereas the neural representation remains independent of the fractional-history mechanism. 
Furthermore, coupled phase-field models contain localized interfacial layers and multiple physical scales, which can lead to an imbalance among the residual, boundary, and data terms during training. Recent phase-field PINN studies indicate that interface-sensitive sampling, normalization, and loss-balancing strategies are important for such coupled systems.
However, the construction of a neural representation that is intrinsically consistent with tempered fractional memory remains largely unexplored.

Motivated by these considerations, we propose \textit{FM-tfPINN}, a physics-informed neural network framework based on fractional memory for tTFCP systems. 
The key idea is to move beyond a purely residual-based treatment of fractional memory by generating coupled phase-field and concentration variables from latent memory-source functions via a tempered fractional-memory mechanism. 
This construction aligns the neural approximation with the history-dependent structure of the governing dynamics before minimizing the residual. 
The resulting framework integrates this memory-generated representation with a fast, shifted residual on graded temporal levels that is accelerated by a SOE, 
as well as interface-aware and residual-adaptive collocation and a unified physics-informed objective for both forward prediction and inverse identification. 
These components are designed to address nonlocal memory, coupled interfacial evolution, concentration redistribution, and sparse physical observations within a single learning framework. 
We examine the performance of the proposed method through representative corrosion-driven phase-field problems, including one-dimensional corrosion-front propagation, distinct kinetic regimes, two-dimensional semi-circular pitting, and inverse mobility identification from physical observations.
These examples are chosen to evaluate FM-tfPINN's ability to recover coupled fields and physically meaningful corrosion diagnostics rather than only pointwise solution values. 
The main features and contributions of this work are: 
\begin{itemize}
\item We develop a fractional-memory 
neural representation for tTFCP systems, where the solution components are constructed through latent memory-source functions.
This formulation incorporates the tempered history effect directly into the neural approximation process.

\item We develop a fast SOE-accelerated shifted residual formulation on graded temporal levels to efficiently enforce the tempered Caputo history contribution in the physics-informed loss.

\item We design interface-aware and residual-adaptive collocation strategies that focus the training process on localized diffuse-interface regions and high-residual zones that arise in corrosion-driven phase-field dynamics.

\item We formulate a unified FM-tfPINN framework for forward and inverse problems by combining coupled residual losses, boundary constraints, admissibility penalties, 
sparse physical observation losses for parameter identification.

\item We demonstrate the framework using physically interpretable diagnostics, including corrosion-depth evolution, interface concentration, concentration peaks, pit geometry measures, tempered mass behavior, and inverse recovery of mobility parameters from sparse observations.
\end{itemize}
The rest of the manuscript is organized as follows. 
Section~\ref{sect:model} introduces the governing tTFCP system, 
with the initial and boundary settings and the associated forward and inverse problem formulations. 
Section~\ref{sect:main_framework} describes the proposed FM-tfPINN methodology, covering the fractional-memory generated neural representation, the fast SOE-accelerated shifted memory residual, the interface-aware and residual-adaptive collocation strategy, and the physics-informed loss formulation for forward prediction and inverse identification. Section~\ref{sect:num_exp} reports numerical experiments on corrosion-driven phase-field dynamics, covering corrosion-front evolution, activation- and diffusion-controlled regimes, semi-circular pitting and inverse mobility identification from sparse physical observations. 
Finally, Section~\ref{sect:concl} concludes the manuscript with a summary of the main findings and a discussion of future research directions.

\section{Mathematical Model and Problem Setting}\label{sect:model}
This section introduces the governing tTFCP system and establishes the notation used in the proposed FM-tfPINN framework. 
The formulation is written in a dimension-independent form, allowing the same model to cover the 1D 
corrosion-front examples, the activation- and diffusion-controlled regimes, the 2D 
semi-circular pitting problem, and the inverse identification setting considered in the numerical section. 
Let $\Omega\subset\mathbb{R}^d$, $d=1,2$, be a bounded spatial domain and let $T>0$ be the final time.
The phase-field variable is denoted by $\phi(\bx,t)$, where $\phi=1$ represents the metal phase and $\phi=0$ represents the electrolyte or pit phase. 
The normalized concentration field is denoted by $c(\bx,t)$. 
For compact notation, we write the coupled state as $\bu(\bx,t)=\bigl(\phi(\bx,t),c(\bx,t)\bigr)^\top$,
$(\bx,t)\in\Omega\times[0,T]$.

\subsection{Governing tTFCP System}\label{subsec:governing_ttfcp_system}
The temporal memory in the tTFCP system is described by the tempered Caputo derivative of order $0<\alpha<1$ with tempering parameter $\lambda\ge0$. 
For a sufficiently smooth scalar function $v$, it is defined by
\begin{equation}\label{eq:model_tempered_caputo}
\partial_t^{\alpha,\lambda}v(\bx,t)= \frac{\mathrm{e}^{-\lambda t}}{\Gamma(1-\alpha)}
\int_0^t (t-s)^{-\alpha}\frac{\partial}{\partial s} 
\bigl(\mathrm{e}^{\lambda s}v(\bx,s)\bigr)\,ds.
\end{equation}
When $\lambda=0$, \eqref{eq:model_tempered_caputo} reduces to the standard Caputo derivative. 
The exponential tempering attenuates the contribution of remote history while retaining the fractional-memory behavior near the initial time. 
The tTFCP system considered in this work is given by
\begin{subequations}\label{eq:model_ttfcp_system}
\begin{align}
\partial_t^{\alpha,\lambda}\phi
&= 2 A_{\rm chem}L_{\rm mob}\bigl[ c-h(\phi)(c_{\rm Se}-c_{\rm Le})-c_{\rm Le}\bigr]
(c_{\rm Se}-c_{\rm Le})h'(\phi)
\nonumber\\
&\qquad -L_{\rm mob}w_{\phi}g'(\phi)+L_{\rm mob}\alpha_{\phi}\Delta\phi,
\label{eq:model_phi}\\
\partial_t^{\alpha,\lambda}c
&=2A_{\rm chem}M_{\rm mob}\Delta c - 2A_{\rm chem}M_{\rm mob}(c_{\rm Se}-c_{\rm Le})\Delta h(\phi).
\label{eq:model_c}
\end{align}
\end{subequations}
for $(\bx,t)\in\Omega\times(0,T]$. 
Here, $A_{\rm chem}$ is the chemical-energy scaling parameter, $L_{\rm mob}$ is the phase-field mobility, $M_{\rm mob}$ is the concentration mobility, $w_{\phi}$ is the double-well coefficient, and $\alpha_{\phi}$ is the gradient-energy coefficient. 
The parameters $c_{\rm Se}$ and $c_{\rm Le}$ represent the equilibrium concentration levels in the solid and liquid phases, respectively. 
The interpolation function and double-well potential are defined by
\begin{equation*}
   h(\phi)=\phi^3(6\phi^2-15\phi+10),\qquad g(\phi)=\phi^2(1-\phi)^2 .
\end{equation*}
The initial conditions are defined as $\phi(\bx,0)=\phi_0(\bx)$,
$c(\bx,0)=c_0(\bx)$,
$\bx\in\Omega$.
The initial functions $\phi_0$ and $c_0$ are prescribed so that they are compatible with the phase-field configuration and the associated concentration distribution. 
In particular, $\phi_0$ represents a diffuse transition between the relevant phases, while $c_0$ is chosen consistently with the interpolation function $h(\phi)$ and the prescribed concentration regime. 
This ensures that the initial state is compatible with the coupled structure of the tTFCP system.
The boundary conditions are represented by
\begin{equation}\label{eq:model_boundary_operator}
\mathcal{B}_{\phi}\phi=0,
\qquad
\mathcal{B}_{c}c=0,
\qquad
\bx\in\partial\Omega,\quad t\in(0,T],
\end{equation}
where $\mathcal{B}_{\phi}$ and $\mathcal{B}_{c}$ denote the boundary operators associated with the phase-field and concentration variables, respectively.
This abstract notation allows the tTFCP system to accommodate different physically relevant boundary settings, including prescribed phase states and homogeneous no-flux constraints. The particular choice of $\mathcal{B}_{\phi}$ and $\mathcal{B}_{c}$ is determined by the physical configuration under consideration, while the FM-tfPINN formulation developed below treats these boundary constraints within a unified physics-informed objective. 
For the development of the FM-tfPINN framework, it is convenient to rewrite \eqref{eq:model_phi}--\eqref{eq:model_c} in the compact form
\begin{equation}\label{eq:model_compact_ttfcp}
\partial_t^{\alpha,\lambda}\bu(\bx,t)
=\bcF[\bu](\bx,t),
\qquad
(\bx,t)\in\Omega\times(0,T],
\end{equation}
where $\bcF[\bu]=\bigl(\mathcal{F}_{\phi}[\phi,c],\mathcal{F}_{c}[\phi,c]\bigr)^\top$
denotes the coupled nonlinear operator defined by the right-hand sides of \eqref{eq:model_phi} and \eqref{eq:model_c}.
This compact form will be used in the next section to construct the fractional-memory generated neural representation and the corresponding physics-informed residual.

\subsection{Forward and inverse problem settings}\label{subsec:forward_inverse_setting}
The forward problem involves determining the coupled fields $\phi$ and $c$ when the physical parameters, initial conditions, boundary conditions, and domain geometry are known. More precisely, for given parameters $\alpha$, $\lambda$, $A_{\rm chem}$, $L_{\rm mob}$, $M_{\rm mob}$, $w_{\phi}$, $\alpha_{\phi}$, $c_{\rm Se}$, and $c_{\rm Le}$, together with the initial data $\phi_0$, $c_0$, and the boundary operators in \eqref{eq:model_boundary_operator}, the forward task is to find the coupled state $\bu=(\phi,c)^\top$ satisfying the tTFCP system \eqref{eq:model_phi}--\eqref{eq:model_c}.
This setting provides the basic prediction problem, assuming all physical parameters are known. The inverse problem is concerned with the identification of unknown physical parameters from limited observations of the system response. 

In this work, the unknown parameter vector is taken as $\bka=(L_{\rm mob},M_{\rm mob})^\top$.
These parameters are selected because they govern the phase-field kinetics and the concentration transport in the tTFCP system. 
Rather than assuming dense measurements of the full solution fields, the inverse setting uses sparse scalar observations associated with physically interpretable corrosion quantities. 
The available observation data are written as
\begin{equation}\label{eq:model_observation_data}
  \mathcal{D}_{\rm obs}=\bigl\{t_m,\,d_{\rm obs}(t_m),\,c_{\Gamma,{\rm obs}}(t_m),\,
  c_{\max,{\rm obs}}(t_m)\bigr\}_{m=1}^{N_{\rm obs}},
\end{equation}
where $d_{\rm obs}$, $c_{\Gamma,{\rm obs}}$, and $c_{\max,{\rm obs}}$ denote the observed corrosion depth, interface concentration, and concentration peak, respectively, and $N_{\rm obs}$ is the number of observation times. 
The inverse task is therefore to recover $\bka$, together with the corresponding coupled fields, while remaining consistent with the governing tTFCP system, the imposed boundary conditions, and the sparse observation data.

This problem setting provides the mathematical foundation for the FM-tfPINN framework, which will be developed in the next section. In the forward setting, the trainable quantities are the neural parameters used to approximate the coupled fields. 
In the inverse setting, the unknown physical parameters in $\bka$ are learned together with the neural parameters by augmenting the physics-informed objective with the observation mismatch associated with $\mathcal{D}_{\rm obs}$. 
Thus, a unified memory-consistent learning formulation allows the same framework to be used for solution prediction and parameter identification.

\section{FM-tfPINN: The Proposed Framework}\label{sect:main_framework}
This section introduces the proposed FM-tfPINN framework for solving the tTFCP system. The construction is designed around four interconnected components.
First, the physical variables are represented by a tempered fractional-memory neural ansatz, in which the coupled fields $(\phi,c)$ are generated from trainable latent memory-source functions. 
Second, the tempered Caputo derivatives are enforced using a fast SOE-accelerated shifted memory residual on a graded temporal grid. 
Third, the residual collocation strategy is adapted to the localized interfacial structure of phase-field corrosion dynamics. 
Lastly, the trainable parameters are determined by a physics-informed objective combining the coupled residual losses, boundary constraints, admissibility penalties, and observation terms for inverse identification, when available.
 Together, these components form a memory-consistent and interface-aware learning framework for forward prediction and parameter recovery in tTFCP systems.
The subsequent subsections detail these components and their integration into the complete FM-tfPINN training procedure.

\subsection{Tempered Fractional-Memory Neural Representation}\label{subsec:fm_generated_representation}
We first introduce the neural representation employed in FM-tfPINN for the tTFCP system \eqref{eq:model_compact_ttfcp}. 
The objective is to incorporate the temporal memory mechanism induced by the tempered Caputo derivative into the trial space, while preserving sufficient flexibility for the neural network to learn the coupled interfacial and concentration dynamics. 
The proposed representation is constructed so that the increments of $\phi$ and $c$ are generated through the same tempered fractional-memory mechanism that appears in the governing tTFCP system. Using the compact form \eqref{eq:model_compact_ttfcp}, the transformed
tempered dynamics suggest that the evolution contains a history contribution
of convolution type. 
In particular, the tempered fractional integral associated with the memory part is defined componentwise by
\begin{equation}\label{eq:tempered_fractional_integral}
\mathcal{I}_t^{\alpha,\lambda}[w](\bx,t)
= \frac{1}{\Gamma(\alpha)}
\int_0^t (t-s)^{\alpha-1}\mathrm{e}^{-\lambda(t-s)} w(\bx,s)\,ds .
\end{equation}
This operator represents the weakly singular and exponentially tempered accumulation of a source term over past times. 
The factor $(t-s)^{\alpha-1}$ characterizes the fractional-memory contribution near the initial time, while $\mathrm{e}^{-\lambda(t-s)}$ controls the attenuation of remote history. 
Hence, the tempering parameter $\lambda$ regulates the effective memory length of the model. 
In the limiting case $\lambda=0$,
\eqref{eq:tempered_fractional_integral} reduces to the standard fractional
integral associated with the Caputo derivative. 
Motivated by this tempered memory operator, the FM-tfPINN trial function is defined by
\begin{equation}\label{eq:fm_tfpinn_coupled_representation}
\bu_{\bth}(\bx,t)
=\bu_b(\bx,t) + \brh(\bx)
\odot\mathcal{I}_t^{\alpha,\lambda}
[\bz_{\bth}](\bx,t),
\end{equation}
where
\begin{equation*}
\bu_{\bth}=(\phi_{\bth},c_{\bth})^\top,
\qquad \bz_{\bth}
   =(z_{\phi,\bth},z_{c,\bth})^\top,
\qquad \brh=(\rho_{\phi},\rho_c)^\top.
\end{equation*}
Here, $\bz_{\bth}$ is a trainable latent neural field, $\bu_b=(\phi_b,c_b)^\top$ is a lifting function selected according to the prescribed initial and boundary data, 
$\brh$ is a componentwise boundary mask, and $\odot$ in \eqref{eq:fm_tfpinn_coupled_representation} denotes the Hadamard product.
Equivalently, \eqref{eq:fm_tfpinn_coupled_representation} gives
\begin{align}\phi_{\bth}(\bx,t)
&=\phi_b(\bx,t)+\rho_{\phi}(\bx)
\mathcal{I}_t^{\alpha,\lambda}[z_{\phi,\bth}](\bx,t),
\label{eq:fm_phi_representation}\\
c_{\bth}(\bx,t)&=c_b(\bx,t)+\rho_c(\bx)
\mathcal{I}_t^{\alpha,\lambda}[z_{c,\bth}](\bx,t).
\label{eq:fm_c_representation}
\end{align}
Thus, the neural network does not directly parameterize the weakly regular physical variables. 
Instead, it learns the latent memory-source fields whose tempered fractional accumulation generates the coupled phase-field and concentration response. 
The masks $\rho_{\phi}$ and $\rho_c$ are chosen according to the boundary structure of the problem. 
The role of the tempering parameter is embedded directly into the representation through the kernel in \eqref{eq:tempered_fractional_integral}.
When $\lambda$ is small, the exponential factor decays slowly and the representation retains a stronger contribution from earlier states. 
When $\lambda$ is larger, the effect of remote history is damped more rapidly. Consequently, the same representation can describe both weakly tempered and strongly tempered corrosion dynamics. 
This is particularly relevant for tTFCP systems, where early-time weak regularity and finite-memory attenuation must be captured simultaneously.

The representation \eqref{eq:fm_tfpinn_coupled_representation} does not rely on prescribing a fixed temporal profile for the solution. 
The weak fractional behavior is induced by the memory operator acting on the learned latent fields. 
Indeed, if $z_{\phi,\bth}$ and $z_{c,\bth}$ are regular near $t=0$, then \eqref{eq:tempered_fractional_integral} gives, for $q\in\{\phi,c\}$,
\begin{equation*} 
\mathcal{I}_t^{\alpha,\lambda}[z_{q,\bth}](\bx,t)
=\frac{z_{q,\bth}(\bx,0)}{\Gamma(1+\alpha)}t^\alpha+\mathcal{O}(t^{\alpha+1}),
\qquad t\to0^+ .
\end{equation*}
Hence,
\begin{equation}\label{eq:learned_coupled_singular_amplitude}
q_{\bth}(\bx,t)-q_b(\bx,t)=\frac{\rho_q(\bx)z_{q,\bth}(\bx,0)}{\Gamma(1+\alpha)}t^\alpha
+\mathcal{O}(t^{\alpha+1}),
\qquad q\in\{\phi,c\}.
\end{equation}
The weak fractional scale is therefore generated by the tempered memory kernel,
while its spatial amplitude and later-time modulation are determined by the trainable latent fields. 
This transfers the approximation task from the weakly regular variables $(\phi,c)$ to the memory-source fields $(z_{\phi,\bth},z_{c,\bth})$, which is more consistent with the structure of tempered fractional evolution. The initial condition is also incorporated naturally. 
Since $\mathcal{I}_0^{\alpha,\lambda}[\bz_{\bth}](\bx,0)=\boldsymbol{0}$,
\eqref{eq:fm_tfpinn_coupled_representation} gives
  $\bu_{\bth}(\bx,0)
  =\bu_b(\bx,0)$.
Thus, choosing $\bu_b(\bx,0)=\bu_0(\bx)$ enforces the prescribed initial state.

\subsection{Fast SOE-Accelerated Shifted Memory Residual}\label{subsec:fast_shifted_memory_residual}
We next describe the discrete memory residual used to enforce the tTFCP system in FM-tfPINN. 
Since the governing equations involve tempered Caputo derivatives, the temporal residual must represent both the weak initial singularity and the nonlocal memory contribution. 
Therefore we employ a graded temporal grid, a shifted high-order approximation of the transformed Caputo history, and an SOE-based compression of the weakly singular kernel. 
Let $0=t_0<t_1<\cdots<t_{N_t}=T$ be a power-law graded mesh: 
\begin{equation}\label{eq:graded_mesh_fm}
  t_n= T\Bigl(\frac{n}{N_t}\Bigr)^r, \qquad 0\le n\le N_t,
\qquad r\ge 1.
\end{equation}
We denote the local temporal step sizes and the maximum temporal step size, respectively, by
$\tau_n=t_n-t_{n-1}$, $1\le n\le N_t$,
with $\tau=\max_{1\le n\le N_t}\tau_n$.
The shifted temporal level is chosen as
\begin{equation*}
  \sigma=1-\frac{\alpha}{2},
  \qquad t_{n+\sigma}=(1-\sigma)t_n+\sigma t_{n+1},
  \qquad 0\le n\le N_t-1 .
\end{equation*}
In forward simulations, $\alpha$ is prescribed and hence $\sigma$ is fixed. In the present work, $\alpha$ is prescribed, and hence the shifted level $\sigma=1-\alpha/2$ is fixed during training. For a scalar quantity $q$, we write $q^n(\bx)=q(\bx,t_n)$.
The tempered Caputo derivative is evaluated through the transformed variable
$\widetilde q(\bx,t)=\mathrm{e}^{\lambda t}q(\bx,t)$.
Accordingly, we define the \textit{tempered increment}
\begin{equation}\label{eq:tempered_difference_fm}
\nabla_\tau^\lambda q^{n+1}
=\widetilde q^{\,n+1}-\widetilde q^{\,n}
=\mathrm{e}^{\lambda t_{n+1}}q^{n+1}-\mathrm{e}^{\lambda t_n}q^n,
\qquad 0\le n\le N_t-1 .
\end{equation}
The shifted nonuniform high-order approximation of the tempered derivative at $t_{n+\sigma}$ can then be written as
\begin{equation}\label{eq:shifted_history_operator_fm}
({}^{\lambda}\mathcal{D}_{\tau}^{\alpha}q)^{n+\sigma}
=\mathrm{e}^{-\lambda t_{n+\sigma}}\sum_{j=0}^{n}
A_{n-j}^{(n+1)}\nabla_\tau^\lambda q^{j+1},
\qquad 0\le n\le N_t-1,
\end{equation}
where the coefficients $A_{n-j}^{(n+1)}$ are generated by the shifted high-order approximation on the graded mesh \cite{DwivediRajeevZeng2026JSC}.
The factor $\mathrm{e}^{-\lambda t_{n+\sigma}}$ maps the Caputo derivative of the transformed variable $\widetilde q$ back to the tempered derivative of $q$.
For weakly singular solutions, the graded mesh \eqref{eq:graded_mesh_fm} is essential for resolving the initial layer.
Under the usual temporal regularity assumptions, such shifted graded approximations yield the consistency behavior
   $\mathcal{O}(\tau^{\min\{r\alpha,2\}})$,
cf.\ \cite{StynesORiordanGracia2017,JinLazarovZhou2016}.

A direct evaluation of \eqref{eq:shifted_history_operator_fm} requires all previous temporal increments and becomes too expensive for long-time simulations and for residual evaluation at many collocation points. 
To avoid this full history summation, we approximate the weakly singular kernel of 
the transformed Caputo derivative by a sum of exponentials \cite{Arnold2003, JiangZhangZhangZhang2017},
\begin{equation}\label{eq:soe_kernel_fm}
\omega_{1-\alpha}(t)
=\frac{t^{-\alpha}}{\Gamma(1-\alpha)}
\approx\sum_{\ell=1}^{N_q}\nu^\ell \mathrm{e}^{-s^\ell t}, \qquad t>0,
\end{equation}
where $\{s^\ell,\nu^\ell\}_{\ell=1}^{N_q}$ are the SOE nodes and weights. 
The tempering is therefore not inserted into the SOE kernel itself. 
Instead, it is incorporated through the transformed increments $\nabla_\tau^\lambda q^{n+1}$ and through the multiplier $\mathrm{e}^{-\lambda t_{n+\sigma}}$ in the final derivative approximation. 
Using \eqref{eq:soe_kernel_fm}, the fast compressed-history approximation is written as
\begin{equation}\label{eq:fast_tempered_operator_fm}
  ({}^{\lambda}\mathcal{D}_{\tau,\mathcal{F}}^{\alpha}q)^{n+\sigma}
=\mathrm{e}^{-\lambda t_{n+\sigma}}
\biggl[a_0^{(n+1)}\nabla_\tau^\lambda q^{n+1}
+\sum_{\ell=1}^{N_q}\nu^\ell\mathrm{e}^{-s^\ell\sigma\tau_{n+1}}
\mathcal{V}_{q}^{\ell}(t_n)\biggr],
\quad
0\le n < N_t.
\end{equation}
Here $\mathcal{V}_{q}^{\ell}$ denotes the compressed history variable associated with the $\ell$-th exponential mode. The memory variables are initialized by $\mathcal{V}_{q}^{\ell}(t_0)=0$, $1\le \ell\le N_q$.
For $1\le m\le N_t-1$, they are updated recursively as
\begin{equation}\label{eq:soe_history_recursion_fm}
   \mathcal{V}_{q}^{\ell}(t_m)=\mathrm{e}^{-s^\ell\tau_m}
   \mathcal{V}_{q}^{\ell}(t_{m-1})
   +a^{(m,\ell)}\nabla_\tau^\lambda q^{m}+b^{(m,\ell)}
   (\varrho_m\nabla_\tau^\lambda q^{m+1}-\nabla_\tau^\lambda q^{m}),
\end{equation}
where $\varrho_m=\tau_m/\tau_{m+1}$.
The coefficients $a^{(m,\ell)}$ and $b^{(m,\ell)}$ are obtained by exact integration of the exponential kernels over the corresponding nonuniform time intervals \cite{JiangZhangZhangZhang2017,DwivediRajeevZeng2026JSC}. 
Through this recursion, the full temporal history is replaced by a finite collection of memory variables.
As a result, only $\mathcal{O}(N_q)$ memory variables need to be stored, and the evaluation of the compressed history contribution at each shifted time level costs $\mathcal{O}(N_q)$ operations.

We now apply \eqref{eq:fast_tempered_operator_fm} to the FM-tfPINN approximation. 
For
$\bu_{\bth}=
(\phi_{\bth},c_{\bth})^\top$,
the shifted neural values are defined componentwise by
\begin{equation*}
q_{\bth}^{n+\sigma}(\bx)
=(1-\sigma)q_{\bth}^{n}(\bx)
+\sigma q_{\bth}^{n+1}(\bx),
\qquad q\in\{\phi,c\}.
\end{equation*}
The corresponding fast tempered derivative of the neural prediction is
\begin{equation}
   ({}^{\lambda}\mathcal{D}_{\tau,\mathcal{F}}^{\alpha}q_{\bth})^{n+\sigma}
=\mathrm{e}^{-\lambda t_{n+\sigma}}
\biggl[ a_0^{(n+1)}\nabla_\tau^\lambda q_{\bth}^{n+1}
+\sum_{\ell=1}^{N_q}\nu^\ell\mathrm{e}^{-s^\ell\sigma\tau_{n+1}}
\mathcal{V}_{q,\bth}^{\ell}(t_n)\biggr],
\quad q\in\{\phi,c\},
\label{eq:fast_tempered_operator_neural_fm}
\end{equation}
where the neural history variables $\mathcal{V}_{q,\bth}^{\ell}$ are generated from the tempered increments of $q_{\bth}$ through the recursion \eqref{eq:soe_history_recursion_fm}. 
Using the compact operator notation introduced in \eqref{eq:model_compact_ttfcp}, 
the componentwise form of the tTFCP system is written as
\begin{equation*}
\partial_t^{\alpha,\lambda}\phi
=\mathcal{F}_{\phi}[\phi,c],
\qquad
\partial_t^{\alpha,\lambda}c
=\mathcal{F}_{c}[\phi,c],
\end{equation*}
where $\mathcal{F}_{\phi}$ and $\mathcal{F}_{c}$ denote the coupled phase-field and concentration operators, respectively. 
Substituting the FM-tfPINN approximation into the shifted compressed-history discretization gives the componentwise residuals
\begin{equation}\label{eq:phi_memory_residual_fm}
  \mathcal{R}_{\phi,\bth}^{n+\sigma}(\bx)
=({}^{\lambda}\mathcal{D}_{\tau,\mathcal{F}}^{\alpha}
\phi_{\bth})^{n+\sigma}(\bx)
-\mathcal{F}_{\phi}[\phi_{\bth}^{n+\sigma},
c_{\bth}^{n+\sigma}](\bx),
\end{equation}
and
\begin{equation}\label{eq:c_memory_residual_fm}
\mathcal{R}_{c,\bth}^{n+\sigma}(\bx)
=({}^{\lambda}\mathcal{D}_{\tau,\mathcal{F}}^{\alpha}
c_{\bth})^{n+\sigma}(\bx)
-\mathcal{F}_{c}\phi_{\bth}^{n+\sigma},
c_{\bth}^{n+\sigma}](\bx),
\end{equation}
for $0\le n\le N_t-1$. The pair
$\boldsymbol{\mathcal{R}}_{\bth}^{n+\sigma}
=(\mathcal{R}_{\phi,\bth}^{n+\sigma},
\mathcal{R}_{c,\bth}^{n+\sigma})^\top$
is referred to as the shifted memory residual.
Its first part contains the local transformed increment and the compressed history contribution, while its second part enforces the coupled phase-field and concentration dynamics at the shifted temporal level.

\subsection{Interface-Aware and Residual-Adaptive Collocation}\label{subsec:interface_adaptive_collocation}
The residuals in \eqref{eq:phi_memory_residual_fm} and \eqref{eq:c_memory_residual_fm} are evaluated on shifted temporal levels and on selected spatial collocation points. 
For the tTFCP system, a uniform distribution of residual points over the whole domain may be inefficient, since the strongest variations of $\phi$ and $c$ are concentrated near the diffuse corrosion interface.
The phase transition, concentration redistribution, and interface motion are all governed by this localized region. 
Therefore, in addition to bulk collocation points, FM-tfPINN employs interface-aware sampling and, for geometrically more complex settings, residual-adaptive enrichment. 
Let $\mathcal{T}_{\sigma}=\{t_{n+\sigma}\colon0\le n\le N_t-1\}$ denote the shifted temporal grid. 
The basic residual collocation set is written as
$\mathcal{S}_{r}=\mathcal{S}_{\rm bulk}\cup\mathcal{S}_{\Gamma_0}$,
where $\mathcal{S}_{\rm bulk}\subset\Omega\times\mathcal{T}_{\sigma}$ contains residual points distributed in the space-time domain, and $\mathcal{S}_{\Gamma_0}$ denotes an enriched set near the initial diffuse interface.
For a one-dimensional corrosion front, the initial interface band can be selected as
\begin{equation*}
  \mathcal{S}_{\Gamma_0}^{1D}=\bigl\{(x,t_{n+\sigma})\colon|x-x_{\Gamma,0}|\le \eta_{\Gamma}\ell,\ 
   t_{n+\sigma}\in\mathcal{T}_{\sigma}\bigr\},
\end{equation*}
where $x_{\Gamma,0}$ is the initial interface location, $\ell$ is the diffuse-interface thickness, and $\eta_{\Gamma}>0$ determines the width of the enriched band. 
For a two-dimensional curved corrosion front, the corresponding initial interface band can be selected according to the initial interface geometry. 
For instance, if the initial front is described by a circular level set with center $(x_c,y_c)$ and radius $R_0$, one may use
\begin{equation*}
  \mathcal{S}_{\Gamma_0}^{2D}=\bigl\{((x,y),t_{n+\sigma})\colon
  \bigl|\sqrt{(x-x_c)^2+(y-y_c)^2}-R_0\bigr|
  \le \eta_{\Gamma}\ell,\ t_{n+\sigma}\in\mathcal{T}_{\sigma}\bigr\}.
\end{equation*}
In corrosion phase-field simulations, the largest spatial gradients and the strongest coupling between $\phi$ and $c$ occur near the diffuse interface.
Enforcing the shifted memory residuals in this region improves the resolution of front propagation, concentration depletion, and interface-induced source terms. 
Thus, the enriched interface band focuses the training process on the part of the domain where the coupled corrosion dynamics are most active.
Boundary collocation sets are introduced separately for boundary conditions that are not imposed directly through the neural representation. 
We denote these sets by $\mathcal{S}_{b,\phi}$ and $\mathcal{S}_{b,c}$. 

Their use depends on the boundary structure of the problem under consideration. 
If a Dirichlet-type constraint is incorporated through a lifting or masking function, the corresponding boundary loss can be omitted.
If no-flux or Neumann-type conditions are prescribed, the required normal-derivative constraints are enforced through boundary collocation points. 
This treatment allows the same FM-tfPINN formulation to accommodate both hard and weak boundary enforcement without modifying the shifted memory residual.

Moreover, the initially enriched set may not remain sufficient during training because the corrosion front evolves and the residual error may become localized away from the initial interface. To address this issue, an additional residual-adaptive enrichment can be used. Given a set of candidate spatial points $\mathcal{C}\subset\Omega$, we compute the residual indicator
\begin{equation}\label{eq:residual_indicator_fm}
\eta_{\bth}(\bx)=\frac{1}{N_t}\sum_{n=0}^{N_t-1}\Bigl(
\bigl|\mathcal{R}_{\phi,\bth}^{n+\sigma}(\bx)\bigr|
+\bigl|\mathcal{R}_{c,\bth}^{n+\sigma}(\bx)\bigr|\Bigr),
\qquad\bx\in\mathcal{C}.
\end{equation}
Points with relatively large values of $\eta_{\bth}$ are selected to form the residual-enriched set $\mathcal{S}_{\rm res}$.
The residual collocation set is then updated as
   $\mathcal{S}_r\leftarrow\mathcal{S}_r \cup \mathcal{S}_{\rm res}$.
This enrichment directs additional residual enforcement to regions where the current approximation has difficulty satisfying the coupled shifted memory equations. 
In geometrically complex corrosion configurations, residual-adaptive enrichment can be complemented by a self-interface enrichment based on the current FM-tfPINN prediction. 
At selected training stages, the predicted diffuse-interface region is identified by
\begin{equation*}
\mathcal{B}_{\Gamma,\bth}(t)
= \{\bx\in\Omega:
\phi_{\min}<\phi_{\bth}(\bx,t)
<\phi_{\max} \},
\end{equation*}
where typical choices are $\phi_{\min}=0.2$ and $\phi_{\max}=0.8$. 
Collocation points sampled from this predicted interface band form the set $\mathcal{S}_{\rm self}$, and the residual set is enriched by
  $\mathcal{S}_{r}\leftarrow \mathcal{S}_{r}
   \cup\mathcal{S}_{\rm self}$.
This procedure tracks the moving diffuse interface using only the current FM-tfPINN approximation.
It does not require reference solutions or exact interface locations. 
The final residual collocation strategy is selected according to the spatial complexity of the corrosion configuration. 
For problems with a 
simple interface motion, the residual set can be taken as
$\mathcal{S}_{r}=\mathcal{S}_{\rm bulk}\cup\mathcal{S}_{\Gamma_0}$.
For problems involving curved or strongly evolving interfaces, the enriched set can be written as
\begin{equation*} 
\mathcal{S}_{r}=\mathcal{S}_{\rm bulk}
\cup\mathcal{S}_{\Gamma_0}\cup\mathcal{S}_{\rm self}\cup\mathcal{S}_{\rm res}.
\end{equation*}
This distinction keeps the collocation strategy flexible: additional enrichment is activated only when the geometry and residual localization require it. To keep the training cost controlled, the enriched residual set is either capped or periodically resampled. In practice, one may impose
 $|\mathcal{S}_{r}|\le N_{r,\max}$,
where $N_{r,\max}$ is a prescribed maximum number of residual points.
This prevents the adaptive enrichment from increasing the computational cost without bound while still allowing the optimizer to focus on the most informative regions.

\subsection{Physics-Informed Loss Formulation}\label{subsec:physics_informed_loss_ttfcp}
We now formulate the training objective used to determine the parameters of the FM-tfPINN approximation. After constructing the tempered fractional-memory representation and the shifted memory residuals, the remaining task is to enforce the coupled phase-field corrosion dynamics, the boundary conditions, and the physically admissible ranges of the predicted variables. 
The loss is therefore assembled from the componentwise residual errors for $\phi$ and $c$, boundary residuals, bound constraints, and, for inverse problems, mismatch terms associated with physical corrosion observations. 
Using the residuals defined in \eqref{eq:phi_memory_residual_fm} and \eqref{eq:c_memory_residual_fm}, we define
\begin{equation}\label{eq:loss_res_phi_fm}
  \mathcal{L}_{r,\phi}= \frac{1}{|\mathcal{S}_r|}
  \sum_{(\bx_i,t_{n+\sigma})\in\mathcal{S}_r}
  \bigl|\mathcal{R}_{\phi,\bth}^{n+\sigma}(\bx_i)\bigr|^2,
\end{equation}
and
\begin{equation}\label{eq:loss_res_c_fm}
   \mathcal{L}_{r,c}=\frac{1}{|\mathcal{S}_r|}
   \sum_{(\bx_i,t_{n+\sigma})\in\mathcal{S}_r}
  \bigl|\mathcal{R}_{c,\bth}^{n+\sigma}(\bx_i)\bigr|^2.
\end{equation}
The total residual loss is written as
    $\mathcal{L}_{r}=\omega_{\phi}\mathcal{L}_{r,\phi}+\omega_{c}\mathcal{L}_{r,c}$,
where $\omega_{\phi}$ and $\omega_c$ balance the contributions of the phase-field and concentration equations. 
This separation is useful because the two equations may have different physical scales, different differential orders, and different sensitivities near the moving diffuse interface. 
We define
\begin{equation}\label{eq:loss_bc_component_fm}
\mathcal{L}_{bc,q}=\frac{1}{|\mathcal{S}_{b,q}|}
\sum_{(\bx_i,t_{n+\sigma})\in\mathcal{S}_{b,q}}
\bigl|\nabla q_{\bth}^{n+\sigma}(\bx_i)
\cdot\boldsymbol{n}-g_{q,b}^{\,n+\sigma}(\bx_i)
\bigr|^2,
\qquad q\in\{\phi,c\},
\end{equation}
where $\boldsymbol{n}$ is the outward unit normal and $g_{q,b}$ denotes the prescribed boundary flux.
For homogeneous no-flux boundaries, $g_{q,b}=0$. 
The total boundary contribution is written as
  $\mathcal{L}_{bc} = \omega_{bc,\phi}\mathcal{L}_{bc,\phi} + \omega_{bc,c}\mathcal{L}_{bc,c}$.
In addition to the governing equations and boundary conditions, 
the corrosion variables must remain physically admissible.
The phase-field variable represents a phase indicator and is expected to remain in the interval $[0,1]$, while the concentration variable should remain nonnegative and, in the nondimensional settings used here, is also restricted to $[0,1]$. 
Therefore we include the bound penalty
\begin{equation}\label{eq:loss_bound_fm}
\begin{aligned}
  \mathcal{L}_{\rm bound}=\frac{1}{|\mathcal{S}_{q}|}
  \sum_{(\bx_i,t_{n+\sigma})\in\mathcal{S}_{q}}\Big[&\ReLU\bigl(-\phi_{\bth}^{n+\sigma}(\bx_i)\bigr)^2
   +\ReLU\bigl(\phi_{\bth}^{n+\sigma}(\bx_i)-1\bigr)^2\\
   &+\ReLU\bigl(-c_{\bth}^{n+\sigma}(\bx_i)\bigr)^2
   +\ReLU\bigl(c_{\bth}^{n+\sigma}(\bx_i)-1\bigr)^2\Big],
\end{aligned}
\end{equation}
where $\mathcal{S}_q$ is the set of points at which the admissibility constraints are monitored.
This term is not intended to replace the physics residual. 
Rather, it stabilizes the optimization by discouraging nonphysical values of $\phi_{\bth}$ and $c_{\bth}$, especially near the diffuse corrosion interface.
For inverse identification from physical corrosion observations, we augment the forward loss by an observation mismatch term. 
Let $\mathcal{S}_{d,\phi}$ and $\mathcal{S}_{d,c}$ denote the available observation sets for the phase-field and concentration variables. 
The field-data losses are defined as
\begin{equation}\label{eq:loss_data_phi_fm}
\mathcal{L}_{d,\phi}
=\frac{1}{|\mathcal{S}_{d,\phi}|}\sum_{(\bx_i,t_i)\in\mathcal{S}_{d,\phi}}
|\phi_{\bth}(\bx_i,t_i)-\phi_i^{\rm obs}|^2,
\end{equation}
and
\begin{equation}\label{eq:loss_data_c_fm}
\mathcal{L}_{d,c}= \frac{1}{|\mathcal{S}_{d,c}|}
\sum_{(\bx_i,t_i)\in\mathcal{S}_{d,c}}|c_{\bth}(\bx_i,t_i)-c_i^{\rm obs}|^2.
\end{equation}
In corrosion applications, observations are often available in the form of geometric or physically aggregated quantities rather than full-field measurements. 
We therefore also allow a diagnostic data term
\begin{equation}\label{eq:loss_data_geo_fm}
   \mathcal{L}_{d,g}=\frac{1}{N_g}\sum_{i=1}^{N_g}|\mathcal{G}_{\bth}(t_i)-\mathcal{G}^{\rm obs}(t_i)|^2,
\end{equation}
where $\mathcal{G}$ may denote corrosion depth, interface concentration, concentration peak, pit area, or equivalent pit radius.
The \textit{total observation loss} is then
\begin{equation}\label{eq:loss_data_total_fm}
\mathcal{L}_{\rm data}
=\omega_{d,\phi}\mathcal{L}_{d,\phi}+\omega_{d,c}\mathcal{L}_{d,c}
+\omega_{d,g}\mathcal{L}_{d,g}.
\end{equation}
For the forward simulations, the data weights are set to zero. 
In the inverse identification setting considered in this work,
the unknown physical parameter vector is taken as
$\bka=(L_{\rm mob},M_{\rm mob})^\top$,
so that the optimization is performed with respect to $(\bth,\bka)$. 
Combining the above contributions, the forward FM-tfPINN objective is
\begin{equation}\label{eq:forward_loss_fm}
   \mathcal{L}_{F}= \mathcal{L}_{r}+ \mathcal{L}_{bc}+ \omega_{\rm bound}\mathcal{L}_{\rm bound}.
\end{equation}
The inverse objective is obtained by adding the observation term,
  $\mathcal{L}_{I}= \mathcal{L}_{F} + \mathcal{L}_{\rm data}$.
Thus, the same loss formulation covers both forward prediction and inverse identification. 
The resulting optimization couples the compressed tempered-memory residuals, interfacial boundary constraints, physical admissibility of the phase and concentration fields, and, when available, physical corrosion observations within a single FM-tfPINN framework.
The complete training workflow is given in Algorithm~\ref{alg:fm_tfpinn}.
All trainable variables, namely the neural parameters $\bth$ and, for inverse problems, the unknown physical parameters $\bka$, are optimized using the Adam algorithm~\cite{KingmaBa2015Adam}. 
The computational implementation is performed in the TensorFlow environment~\cite{AbadiEtAl2016TensorFlow}. 

\begin{algorithm}
\scriptsize
\caption{FM-tfPINN training procedure for the tTFCP system}
\label{alg:fm_tfpinn}
\begin{algorithmic}[1]
\Require Domain $\Omega$, final time $T$, time levels $N_t$, grading exponent $r$, parameters $\alpha,\lambda$, initial state $\bu_0=(\phi_0,c_0)^\top$, lifting $\bu_b$, masks $\brh$, collocation sets $\mathcal{S}_{\rm bulk},\mathcal{S}_{\Gamma_0},\mathcal{S}_{b,\phi},\mathcal{S}_{b,c},\mathcal{S}_{q}$, optional adaptive sets $\mathcal{S}_{\rm self},\mathcal{S}_{\rm res}$, candidate set $\mathcal{C}$, maximum residual-set size $N_{r,\max}$, optional observations $\mathcal{D}_{\rm obs}$, SOE nodes and weights $\{(s^\ell,\nu^\ell)\}_{\ell=1}^{N_q}$, loss weights, maximum iterations $M_{\mathrm{iter}}$, tolerance $\varepsilon_{\mathrm{tol}}$, learning rate $\eta_{\mathrm{lr}}$
\Ensure Trained parameters $\bth^{*}$, optional inverse parameters $\bka^{*}$, and learned approximation $\bu_{\bth^{*}}=(\phi_{\bth^{*}},c_{\bth^{*}})^\top$

\State Initialize neural parameters $\bth^{(0)}$ for $\bz_{\bth}=(z_{\phi,\bth},z_{c,\bth})^\top$
\If{inverse identification is considered}
    \State Initialize $\bka^{(0)}=(L_{\rm mob}^{(0)},M_{\rm mob}^{(0)})^\top$
\EndIf
\State Build $t_n=T(n/N_t)^r$, $0\le n\le N_t$, and set $\sigma=1-\alpha/2$, $t_{n+\sigma}=(1-\sigma)t_n+\sigma t_{n+1}$
\State Generate shifted temporal coefficients and SOE quadrature data
\State Set the initial residual set $\mathcal{S}_{r}=\mathcal{S}_{\rm bulk}\cup\mathcal{S}_{\Gamma_0}$

\For{$it=1,\ldots,M_{\mathrm{iter}}$}
    \If{inverse identification is considered}
    \State Update the physical coefficients in the residual using the current $\bka$
    \EndIf

    \If{self-interface enrichment is activated at iteration $it$}
        \State Sample $\mathcal{S}_{\rm self}$ from $\{\bx\in\Omega:\phi_{\min}<\phi_{\bth}(\bx,t)<\phi_{\max}\}$ and set $\mathcal{S}_r\leftarrow\mathcal{S}_r\cup\mathcal{S}_{\rm self}$ for subsequent residual evaluations
    \EndIf
    \If{residual-adaptive enrichment is activated at iteration $it$}
        \State Select $\mathcal{S}_{\rm res}$ from $\mathcal{C}$ using $\eta_{\bth}(\bx)$ in \eqref{eq:residual_indicator_fm}, set $\mathcal{S}_r\leftarrow\mathcal{S}_r\cup\mathcal{S}_{\rm res}$, and cap $|\mathcal{S}_r|\le N_{r,\max}$ for subsequent residual evaluations
    \EndIf

    \State Evaluate $\bu_{\bth}=\bu_b+\brh\odot\mathcal{I}_t^{\alpha,\lambda}[\bz_{\bth}]$ on the active residual, boundary, admissibility, and data sets
    \State Form $q_{\bth}^{n+\sigma}=(1-\sigma)q_{\bth}^{n}+\sigma q_{\bth}^{n+1}$, $q\in\{\phi,c\}$
    \State Compute $\nabla_\tau^\lambda q_{\bth}^{n+1}=\mathrm{e}^{\lambda t_{n+1}}q_{\bth}^{n+1}-\mathrm{e}^{\lambda t_n}q_{\bth}^{n}$
    \State Reinitialize $\mathcal{V}^{\ell}_{q,\bth}(t_0)=0$, $q\in\{\phi,c\}$, $1\le \ell\le N_q$, for the current network parameters
    \State Recompute the SOE histories $\mathcal{V}^{\ell}_{q,\bth}$ from the current tempered increments using \eqref{eq:soe_history_recursion_fm}
    \State Assemble $({}^{\lambda}\mathcal{D}^{\alpha}_{\tau,\mathcal{F}}q_{\bth})^{n+\sigma}$ using \eqref{eq:fast_tempered_operator_neural_fm}
    \State Compute $\mathcal{R}_{\phi,\bth}^{n+\sigma}$ and $\mathcal{R}_{c,\bth}^{n+\sigma}$ from \eqref{eq:phi_memory_residual_fm}--\eqref{eq:c_memory_residual_fm}
    \State Evaluate $\mathcal{L}_{r}$, $\mathcal{L}_{bc}$, $\mathcal{L}_{\rm bound}$, and, if observations are available, $\mathcal{L}_{\rm data}$
    \State Set $\mathcal{L}=\mathcal{L}_{F}$ for forward prediction and $\mathcal{L}=\mathcal{L}_{I}$ for inverse identification

    \If{inverse identification is considered}
        \State Update $(\bth,\bka)\leftarrow(\bth,\bka)-\eta_{\mathrm{lr}}\nabla_{(\bth,\bka)}\mathcal{L}$
    \Else
        \State Update $\bth\leftarrow\bth-\eta_{\mathrm{lr}}\nabla_{\bth}\mathcal{L}$
    \EndIf

    \If{$\mathcal{L}<\varepsilon_{\mathrm{tol}}$}
        \State \textbf{break}
    \EndIf
\EndFor

\State Set $\bth^{*}=\bth$ and, if applicable, $\bka^{*}=\bka$
\State \Return $\bth^{*}$, $\bka^{*}$ if applicable, and $\bu_{\bth^{*}}$
\end{algorithmic}
\end{algorithm}

\section{Numerical Experiments}\label{sect:num_exp}
Here we assess the performance of the proposed FM-tfPINN framework for forward prediction and inverse identification of the tTFCP system.
The numerical tests are designed to examine the ability of the method to capture weak initial fractional regularity, tempered memory attenuation, moving diffuse interfaces, concentration-driven corrosion response and geometry-dependent pitting evolution. 
All experiments employ the graded temporal mesh $t_n=T(n/N_t)^r$, the shifted levels $t_{n+\sigma}$, the SOE-compressed tempered memory residual and the fractional-memory generated approximation $\bu_{\bth}=(\phi_{\bth},c_{\bth})^\top$. 
The network architecture, collocation size, loss weights, training iterations, and example-specific physical parameters are reported separately for each test. 
The error measures and physical diagnostics used throughout the numerical section are summarized in Table~\ref{tab:numerical_diagnostics}.

\begin{table}[!t]
\centering
\scriptsize
\caption{Error measures and physical diagnostics used in the numerical experiments.}
\label{tab:numerical_diagnostics}
\small
\renewcommand{\arraystretch}{1.25}
\begin{tabularx}{\textwidth}{@{}p{0.24\textwidth}X p{0.25\textwidth}@{}}
\toprule
\textbf{Category} & \textbf{Quantity} & \textbf{Used in} \\
\midrule

Field accuracy
& $E_2^q=\frac{\|q_{\bth}-q_{\rm ref}\|_2}
{\|q_{\rm ref}\|_2},\,E_\infty^q=\|q_{\bth}-q_{\rm ref}\|_\infty,$ and $q\in\{\phi,c\}$.
&All forward examples\\

One-dimensional corrosion response
&Corrosion depth $d(t)$, interface concentration $c_{\Gamma}(t)$, and concentration peak $c_{\max}(t)$.
&Examples \ref{subsec:ex1_pencil_electrode} and \ref{subsec:ex2_activation_diffusion}\\

Tempered-memory behavior
& Ordinary concentration mass $M_c(t)$ and tempered mass $\mathrm{e}^{\lambda t}M_c(t)$.
& Memory-consistency check\\

Two-dimensional pitting response
& Pit area $A_{\rm pit}(t)$ and equivalent radius $R_{\rm eq}(t)$.
& Example \ref{subsec:ex3_semicircular_pitting}\\

Inverse identification
& Absolute parameter error $|\kappa_j-\kappa_{j,\rm ref}|$, relative parameter error
$|\kappa_j-\kappa_{j,\rm ref}|/|\kappa_{j,\rm ref}|$, and observation mismatch $E_{\rm obs}$.
& Example~\ref{subsec:ex4_inverse_identification}\\

\bottomrule
\end{tabularx}
\end{table}

\subsection{Tempered fractional Pencil-Electrode Corrosion}\label{subsec:ex1_pencil_electrode}

\subsubsection{Problem Statement}\label{subsubsec:ex1_problem_statement}
We first consider a one-dimensional tempered fractional pencil-electrode corrosion problem on
$\Omega=(-H_s,H_l)$.
The two coupled unknowns are the phase-field variable $\phi(x,t)$ and the normalized concentration field $c(x,t)$. 
In this test, $\phi=1$ corresponds to the metal region, whereas $\phi=0$ corresponds to the electrolyte region. 
For $0<\alpha<1$, $\lambda\ge0$, the governing tTFCP system  reads
\begin{equation}\label{eq:ex1_pencil_model}
\begin{cases}
\begin{aligned}
\partial_t^{\alpha,\lambda}\phi
={}&2A_{\rm chem}L_{\rm mob}\bigl[c-h(\phi)(c_{\rm Se}-c_{\rm Le})-c_{\rm Le}\bigr]
(c_{\rm Se}-c_{\rm Le})h'(\phi)\\
& -L_{\rm mob}w_{\phi}g'(\phi)+L_{\rm mob}\alpha_{\phi}\phi_{xx},
\end{aligned}\\
\begin{aligned}
\partial_t^{\alpha,\lambda}c={}& 2A_{\rm chem}M_{\rm mob}c_{xx}
  -2A_{\rm chem}M_{\rm mob}(c_{\rm Se}-c_{\rm Le})(h(\phi))_{xx},\\
&\qquad x\in(-H_s,H_l),\quad t\in(0,T],
\end{aligned}
\\
\phi(x,0)=\phi_0(x),
\qquad c(x,0)=c_0(x),
\qquad x\in[-H_s,H_l],\\
\phi(-H_s,t)=1,
\qquad \phi(H_l,t)=0,
\qquad t\in[0,T],
\\ c_x(-H_s,t)=0,
\qquad c_x(H_l,t)=0,
\qquad t\in[0,T].
\end{cases}
\end{equation}
Here, $A_{\rm chem}$ denotes the chemical-energy scaling factor, while $L_{\rm mob}$ and $M_{\rm mob}$ specify the mobilities associated with the phase-field and concentration equations, respectively. 
The parameter $w_{\phi}$ controls the double-well energy contribution, and $\alpha_{\phi}$ denotes the gradient-energy coefficient.
The quantities $c_{\rm Se}$ and $c_{\rm Le}$ represent the equilibrium concentration levels corresponding to the solid and liquid phases, respectively. 
The interpolation function and double-well potential are chosen as
\begin{equation*}
  h(\phi)=\phi^3(6\phi^2-15\phi+10),\qquad g(\phi)=\phi^2(1-\phi)^2 .
\end{equation*}
The initial diffuse interface is centered at $x=0$ and is prescribed by
\begin{equation}\label{eq:ex1_initial_condition}
  \phi_0(x)=\frac{1}{2}\biggl[1-\tanh\Bigl(\frac{x}{\sqrt{2}\,\ell}\Bigr)\biggr],
\qquad c_0(x)=h(\phi_0(x))c_{\rm Se},
\qquad \ell=\sqrt{\frac{\alpha_{\phi}}{w_{\phi}}}.
\end{equation}
The Dirichlet boundary conditions for $\phi$ fix the metal and electrolyte states at the two ends of the computational interval, while the homogeneous Neumann conditions for $c$ impose zero concentration flux through the boundary.
This example provides a physical forward test for FM-tfPINN, since the method must resolve a moving diffuse corrosion front, concentration redistribution, weak fractional-memory behavior near the initial time, and the finite-memory attenuation induced by the tempering factor.

\subsubsection{Implementation Details}
\label{subsubsec:ex1_implementation_details}
The numerical treatment of this problem is carried out using the FM-tfPINN procedure described in Algorithm~\ref{alg:fm_tfpinn}. 
In accordance with the fractional-memory construction introduced in Section~\ref{subsec:fm_generated_representation}, the neural approximation is written in the form
\begin{equation}\label{eq:ex1_fm_representation}
\phi_{\bth}(x,t)=\phi_0(x)
+\rho_{\phi}(x)\mathcal{I}_t^{\alpha,\lambda}[z_{\phi,\bth}](x,t),
\quad c_{\bth}(x,t)
= c_0(x)+\mathcal{I}_t^{\alpha,\lambda}[z_{c,\bth}](x,t),
\end{equation}
where the spatial mask for the phase-field component is chosen as 
\begin{equation*} 
   \rho_{\phi}(x)=(x+H_s)(H_l-x).
\end{equation*}
This construction incorporates the initial profiles directly, since $\mathcal{I}_t^{\alpha,\lambda}[\cdot]=0$ at $t=0$.
In addition, the factor $\rho_{\phi}$ vanishes at both endpoints of the computational interval, so the boundary values of the phase-field lifting are retained throughout training. 
For the present choice of $\phi_0$, these endpoint values are consistent with the prescribed metal and electrolyte states up to negligible truncation error. The concentration equation is treated with the same memory-generated form, while its homogeneous no-flux boundary condition is imposed through the boundary loss. 
Since the present test is a forward prediction problem, no observation term is included in the objective.
A reference solution is computed with a fast SOE-accelerated shifted finite-difference discretization of the same tempered fractional coupled phase-field model. 
The FM-tfPINN residual is evaluated at the shifted temporal levels $\{t_{n+\sigma}\}_{n=0}^{N_t-1}$, where $\sigma=1-\alpha/2$, and at spatial collocation points consisting of bulk samples together with additional points concentrated near the initial diffuse interface. The main computational parameters used in this experiment are given in Table~\ref{tab:ex1_implementation_details}.

\begin{table}[!t]
\centering
\caption{Implementation details for Example~\ref{subsec:ex1_pencil_electrode}.}
\label{tab:ex1_implementation_details}
\scriptsize
\renewcommand{\arraystretch}{1.18}
\begin{tabular}{ll}
\toprule
\textbf{Item} & \textbf{Setting} \\
\midrule
Configuration & FM-tfPINN forward prediction \\
Computational domain & $\Omega=(-1,1)$, with $H_s=H_l=1$ \\
Final time & $T=0.2$ \\
Fractional order and tempering & $\alpha=0.8,\quad \lambda=1$ \\
Shifted level & $\sigma=1-\alpha/2=0.6$ \\
FM-tfPINN temporal grid & $N_t=32,\quad r=(3-\alpha)/\alpha$ \\
Reference temporal grid & $N_t^{\rm ref}=800,\quad r_{\rm ref}=2$ \\
Memory-integral quadrature & $8$-point Gauss--Legendre rule \\
Residual spatial points & $30$ bulk points and $30$ initial-interface enriched points \\
Interface band width & $3\ell$, where $\ell=\sqrt{\alpha_{\phi}/w_{\phi}}$ \\
Neural network & $3$ hidden layers with $15$ neurons per layer \\
Network input and output & $(x,t)\mapsto(z_{\phi,\bth},z_{c,\bth})$ \\
Activation function & $\tanh$ \\
Physical parameters &
$A_{\rm chem}=1,\ L_{\rm mob}=M_{\rm mob}=10^{-3},\
w_{\phi}=1,\ \alpha_{\phi}=10^{-3}$ \\
Equilibrium concentrations & $c_{\rm Se}=1,\quad c_{\rm Le}=0.036$ \\
Optimizer & Adam \\
Learning rate & $10^{-4}$ \\
Training iterations & $3.0\times10^{4}$ \\
Loss weights &
$\omega_{\phi}=1,\ \omega_c=1,\ \omega_{bc,c}=1,\
\omega_{\rm bound}=5\times10^{-2},\ \omega_d=0$ \\
\bottomrule
\end{tabular}
\renewcommand{\arraystretch}{1}
\end{table}

\subsubsection{Results and Discussion}\label{subsubsec:ex1_results_discussion}
\begin{figure}[htbp]
  \centering
  \begin{subfigure}[b]{0.64\textwidth}
    \includegraphics[width=\textwidth]{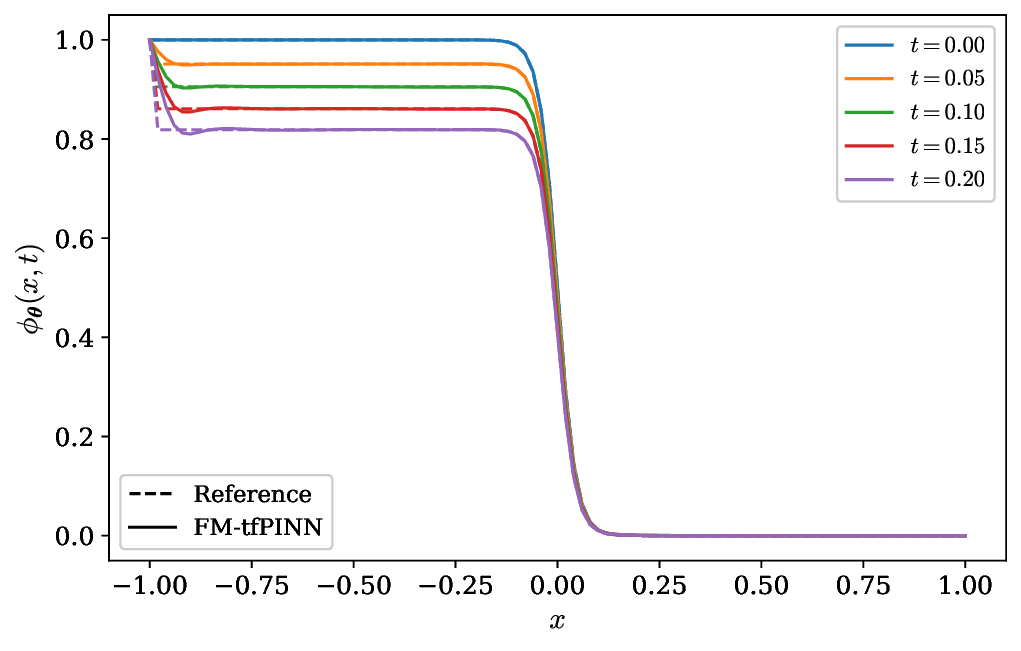}
    \caption{}
    \label{fig:ex1_phi_profiles}
  \end{subfigure}
  \hfill
  \begin{subfigure}[b]{0.64\textwidth}
    \includegraphics[width=\textwidth]{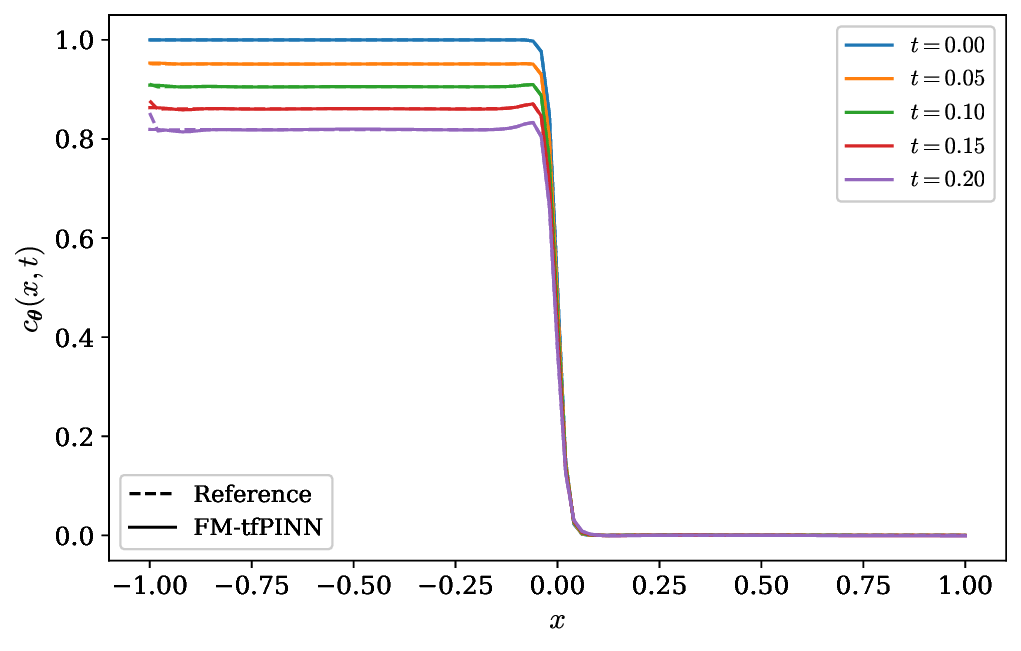}
    \caption{}
    \label{fig:ex1_c_profiles}
  \end{subfigure}
  \caption{Profile comparison for Example~\ref{subsec:ex1_pencil_electrode}. 
  Subfigure (a) shows the phase-field profile $\phi_{\bth}(x,t)$ against the reference solution. 
  Subfigure (b) shows the concentration profile $c_{\bth}(x,t)$ against the reference solution.}
  \label{fig:ex1_profiles}
\end{figure}

We now examine the predictive accuracy and physical consistency of FM-tfPINN for the tempered fractional pencil-electrode corrosion problem. 
Figure~\ref{fig:ex1_profiles} compares the learned profiles with the reference solution at representative time levels.
The predicted phase-field $\phi_{\bth}$ accurately follows the reference interfacial transition,
while the concentration field $c_{\bth}$ captures the coupled redistribution induced by the moving corrosion front. 
The agreement is maintained in the interfacial region, where the spatial gradients are strongest and where the coupling between $\phi$ and $c$ is most pronounced. 
This confirms that the fractional-memory neural representation can resolve the coupled tempered dynamics without using observation data in the forward setting.

The quantitative results are summarized in Table~\ref{tab:ex1_diagnostics}. 
The relative $L^2$ errors for both components are of order $10^{-3}$, which indicates that the proposed method gives an accurate global approximation of the coupled solution. 
The maximum errors remain localized and are mainly associated with the diffuse-interface transition, where small spatial shifts in the interface location can produce larger pointwise deviations.
The physical corrosion-depth error is also small, showing that FM-tfPINN captures not only the field variables but also the derived moving-interface quantity.
In addition, the tempered concentration-mass drift remains below $10^{-3}$, providing a memory-consistency check for the learned solution.

\begin{table}[htbp]
\centering
\scriptsize
\caption{Quantitative diagnostics for Example~\ref{subsec:ex1_pencil_electrode}.}
\label{tab:ex1_diagnostics}
\renewcommand{\arraystretch}{1.25}
\begin{tabular}{lc}
\toprule
Diagnostic quantity & Value \\
\midrule
Relative $L^2$ error of the phase field, $E_2^\phi$ & $1.960235\mathrm{e}{-}03$ \\
Relative $L^2$ error of the concentration field, $E_2^c$ & $1.997110\mathrm{e}{-}03$ \\
Maximum error of the phase field, $E_\infty^\phi$ & $1.065786\mathrm{e}{-}02$ \\
Maximum error of the concentration field, $E_\infty^c$ & $1.195869\mathrm{e}{-}02$ \\
Relative error of the corrosion depth, $d(t)$ & $5.597786\mathrm{e}{-}03$ \\
Predicted final corrosion depth, $d_{\bth}(t_{N_t-1+\sigma})$ & $1.052768\mathrm{e}{-}02$ \\
Reference final corrosion depth, $d_{\rm ref}(t_{N_t-1+\sigma})$ & $1.058363\mathrm{e}{-}02$ \\
Maximum tempered concentration-mass drift & $4.798906\mathrm{e}{-}04$ \\
\bottomrule
\end{tabular}
\renewcommand{\arraystretch}{1}
\end{table}

\begin{figure}[H]
  \centering
  \begin{subfigure}[b]{0.48\textwidth}
    \includegraphics[width=\textwidth]{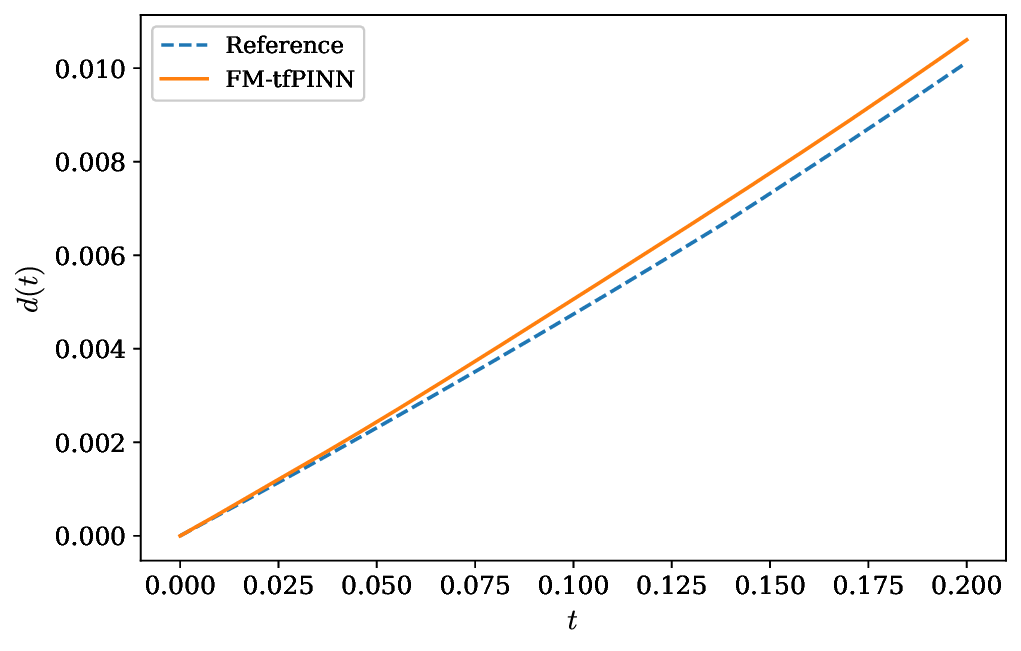}
    \caption{}
    \label{fig:ex1_corrosion_depth}
  \end{subfigure}
  \hfill
  \begin{subfigure}[b]{0.48\textwidth}
    \includegraphics[width=\textwidth]{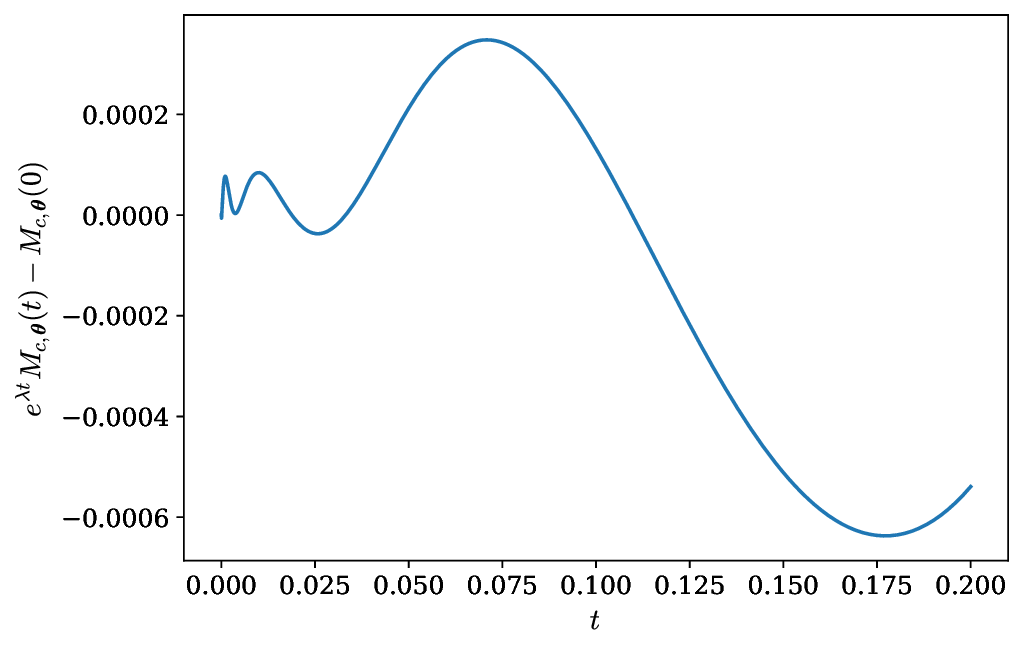}
    \caption{}
    \label{fig:ex1_tempered_mass_drift}
  \end{subfigure}
  \caption{Physical diagnostics for Example~\ref{subsec:ex1_pencil_electrode}. Subfigure (a) compares the predicted corrosion depth $d_{\bth}(t)$ with the reference depth $d_{\rm ref}(t)$. Subfigure (b) shows the tempered concentration-mass drift $\mathrm{e}^{\lambda t}M_{c,\bth}(t)-M_{c,\bth}(0)$.}
  \label{fig:ex1_depth_mass}
\end{figure}

Figure~\ref{fig:ex1_depth_mass} further reports the physical diagnostics associated with corrosion-front motion and tempered mass behavior. The predicted corrosion depth $d_{\bth}(t)$ closely follows the reference depth $d_{\rm ref}(t)$ throughout the computational time interval.
This is a nontrivial validation because the depth is extracted from the learned diffuse interface and is not imposed directly during training. 
The result therefore demonstrates that the learned solution preserves the physically relevant interface motion generated by the coupled tTFCP dynamics.

The tempered concentration-mass diagnostic in Figure~\ref{fig:ex1_depth_mass} provides an additional check on the memory-aware structure of the approximation. 
The small drift of $\mathrm{e}^{\lambda t}M_{c,\bth}(t)-M_{c,\bth}(0)$ indicates that the learned concentration field remains consistent with the expected tempered-memory response. 
Overall, this example shows that FM-tfPINN can simultaneously approximate the coupled phase-field and concentration variables, recover corrosion-front motion and maintain a stable tempered-memory diagnostic for a one-dimensional pencil-electrode corrosion process.

\subsection{Activation- and Diffusion-Controlled Corrosion Regimes}\label{subsec:ex2_activation_diffusion}

\subsubsection{Problem Statement}\label{subsubsec:ex2_problem_statement}
We next consider a regime-dependent corrosion test governed by the same tTFCP system introduced in Example~\ref{subsec:ex1_pencil_electrode}. The purpose of this experiment is to examine whether FM-tfPINN can resolve two qualitatively different corrosion mechanisms within the same coupled phase-field formulation. 
The first case represents an activation-controlled regime, where the phase-field mobility is small and the interfacial reaction is comparatively slow. 
The second case represents a diffusion-controlled regime, where the phase-field mobility is increased so that the concentration transport becomes the dominant limiting mechanism. 
Thus, the model equations, interpolation function $h(\phi)$, double-well potential $g(\phi)$, initial diffuse interface, and boundary conditions are kept identical to those in \eqref{eq:ex1_pencil_model}--\eqref{eq:ex1_initial_condition}, while the corrosion mechanism is distinguished through the mobility pair
\begin{equation}\label{eq:ex2_regime_mobility}
(L_{\rm mob},M_{\rm mob})=\begin{cases}
     (10^{-4},10^{-3}), & \text{activation-controlled regime},\\
     (5\times10^{-1},10^{-3}), & \text{diffusion-controlled regime}.
\end{cases}
\end{equation}
In both cases, the concentration mobility $M_{\rm mob}$ is fixed, while the phase-field mobility $L_{\rm mob}$ is varied to produce distinct corrosion responses.
This configuration provides a more demanding validation than a single-regime forward test,
since the method must capture not only the coupled evolution of $\phi$ and $c$, but also the changes in corrosion depth, interface concentration, and concentration peak induced by the underlying kinetic regime.

\subsubsection{Implementation Details}\label{subsubsec:ex2_implementation_details}
The numerical implementation follows the FM-tfPINN procedure in Algorithm~\ref{alg:fm_tfpinn}. Since the governing equations, boundary conditions, and initial profiles are inherited from Example~\ref{subsec:ex1_pencil_electrode}, the same fractional-memory neural representation in \eqref{eq:ex1_fm_representation} is used for both regimes.
In particular, the phase-field component is represented with the endpoint-vanishing mask $\rho_\phi(x)$, and the concentration component is represented through the tempered fractional-memory increment. 
The Dirichlet boundary conditions for $\phi$ are therefore imposed while the homogeneous no-flux condition for $c$ is imposed through the boundary loss. 
No observation data are used in this forward experiment. 

For each regime, a corresponding reference solution is generated using the same fast SOE-accelerated shifted finite-difference solver used in Example~\ref{subsec:ex1_pencil_electrode}. 
The FM-tfPINN residual is evaluated on the shifted temporal levels $t_{n+\sigma}$, with $\sigma=1-\alpha/2$, and the residual collocation set includes both bulk points and points concentrated near the initial diffuse interface. 
In addition to the field errors $E_2^\phi$, $E_2^c$, $E_\infty^\phi$, and $E_\infty^c$, this example reports the corrosion depth $d(t)$, the interface concentration $c_\Gamma(t)$, and the concentration peak $c_{\max}(t)$, which are the relevant physical diagnostics for distinguishing the two regimes. 
The settings are summarized in Table~\ref{tab:ex2_implementation_details}.

\begin{table}[htbp]
\centering
\scriptsize
\caption{Regime-dependent implementation details for Example~\ref{subsec:ex2_activation_diffusion}. 
The remaining model functions, initial profiles, boundary conditions, optimizer, and neural representation are the same as in Example~\ref{subsec:ex1_pencil_electrode}.}
\label{tab:ex2_implementation_details}
\renewcommand{\arraystretch}{1.25}
\begin{tabular}{lll}
\toprule
Item & Activation-controlled regime & Diffusion-controlled regime \\
\midrule
Configuration & FM-tfPINN forward prediction & FM-tfPINN forward prediction \\
Tempering parameter & $\lambda=1$ & $\lambda=1$ \\
Fractional order & $\alpha=0.8$ & $\alpha=0.8$ \\
Phase-field mobility & $L_{\rm mob}=10^{-4}$ & $L_{\rm mob}=5\times10^{-1}$ \\
Concentration mobility & $M_{\rm mob}=10^{-3}$ & $M_{\rm mob}=10^{-3}$ \\
Shifted level & $\sigma=1-\alpha/2=0.6$ & $\sigma=1-\alpha/2=0.6$ \\
FM-tfPINN temporal grid & $N_t=20,\quad r=(3-\alpha)/\alpha$ & $N_t=20,\quad r=(3-\alpha)/\alpha$ \\
Main diagnostics & $d(t),\ c_\Gamma(t),\ c_{\max}(t)$ & $d(t),\ c_\Gamma(t),\ c_{\max}(t)$ \\
Training iterations & $3.0\times10^{4}$ & $3.0\times10^{4}$ \\
Data weight & $\omega_d=0$ & $\omega_d=0$ \\
\bottomrule
\end{tabular}
\renewcommand{\arraystretch}{1}
\end{table}

\subsubsection{Results and Discussion}\label{subsubsec:ex2_results_discussion}
\begin{figure}[htbp]
  \centering
  \begin{subfigure}[b]{0.48\textwidth}
    \includegraphics[width=\textwidth]{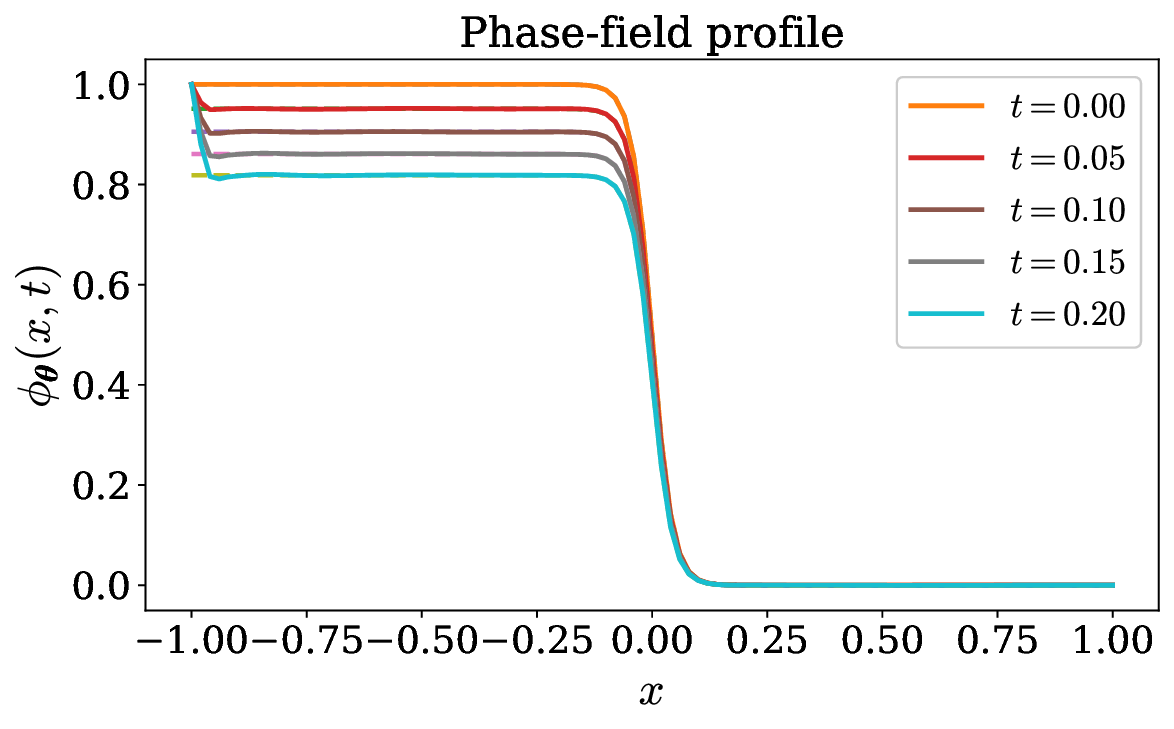}
    \caption{}
    \label{fig:ex2_activation_phi_profiles}
  \end{subfigure}
  \hfill
  \begin{subfigure}[b]{0.48\textwidth}
    \includegraphics[width=\textwidth]{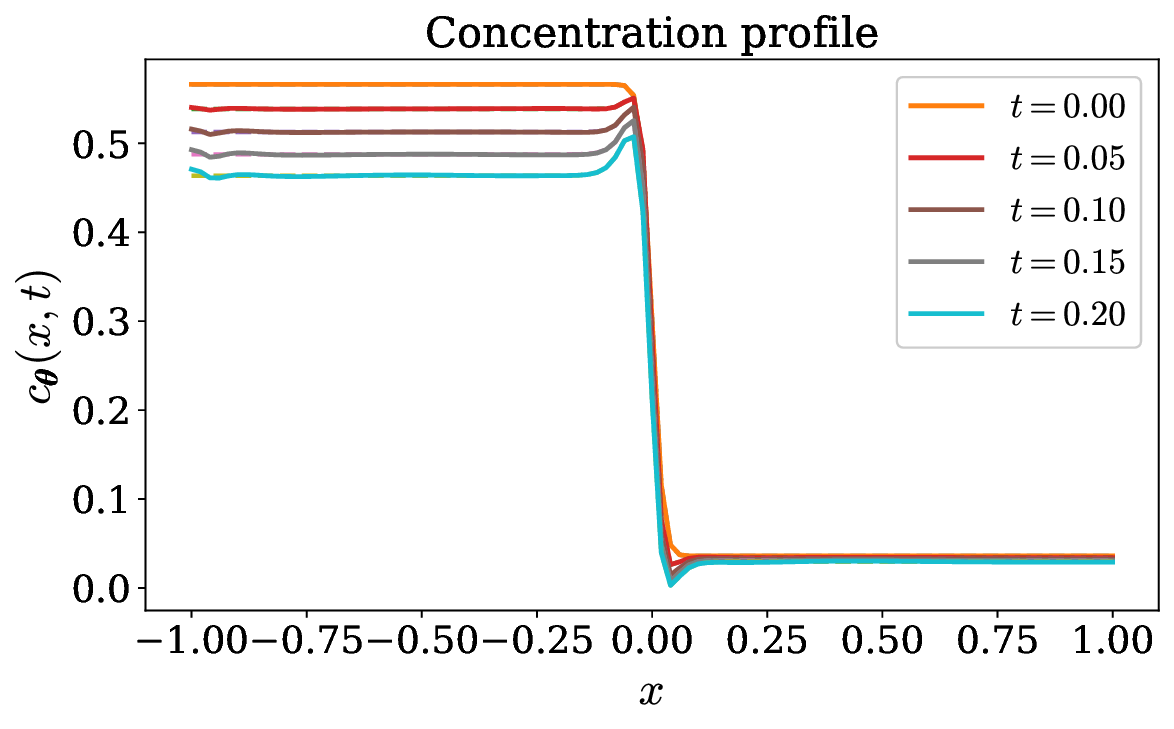}
    \caption{}
    \label{fig:ex2_activation_c_profiles}
  \end{subfigure}

  \begin{subfigure}[b]{0.48\textwidth}
    \includegraphics[width=\textwidth]{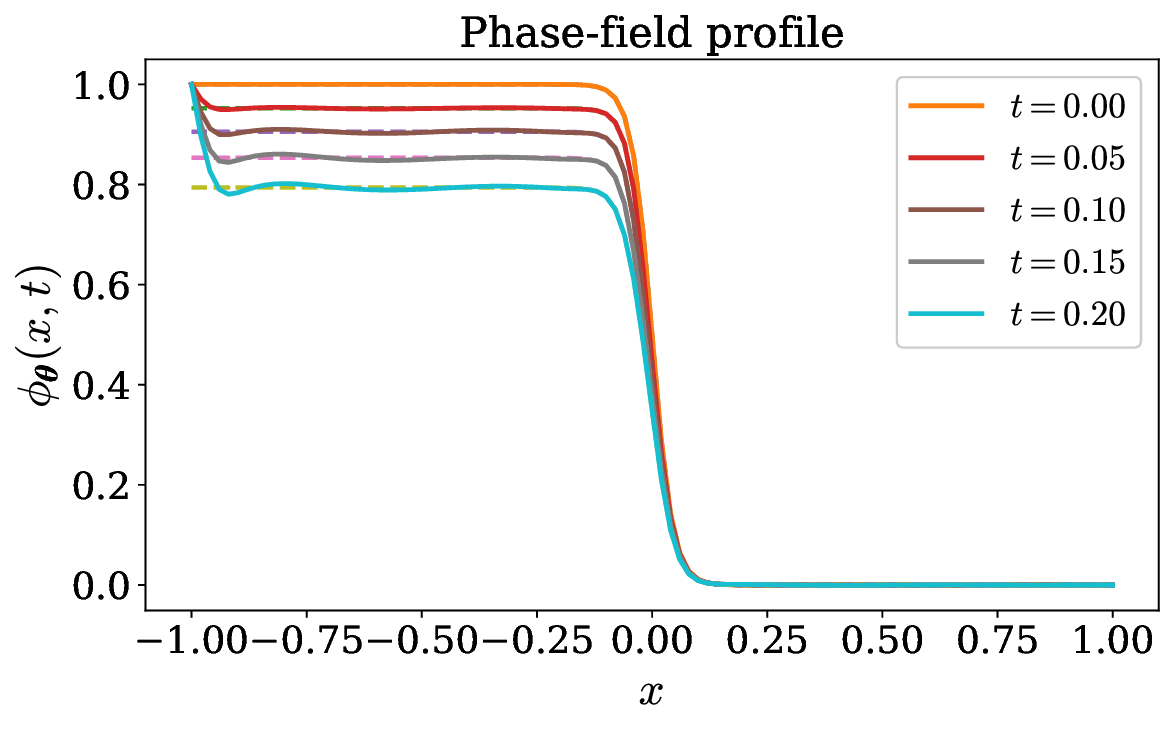}
    \caption{}
    \label{fig:ex2_diffusion_phi_profiles}
  \end{subfigure}
  \hfill
  \begin{subfigure}[b]{0.48\textwidth}
    \includegraphics[width=\textwidth]{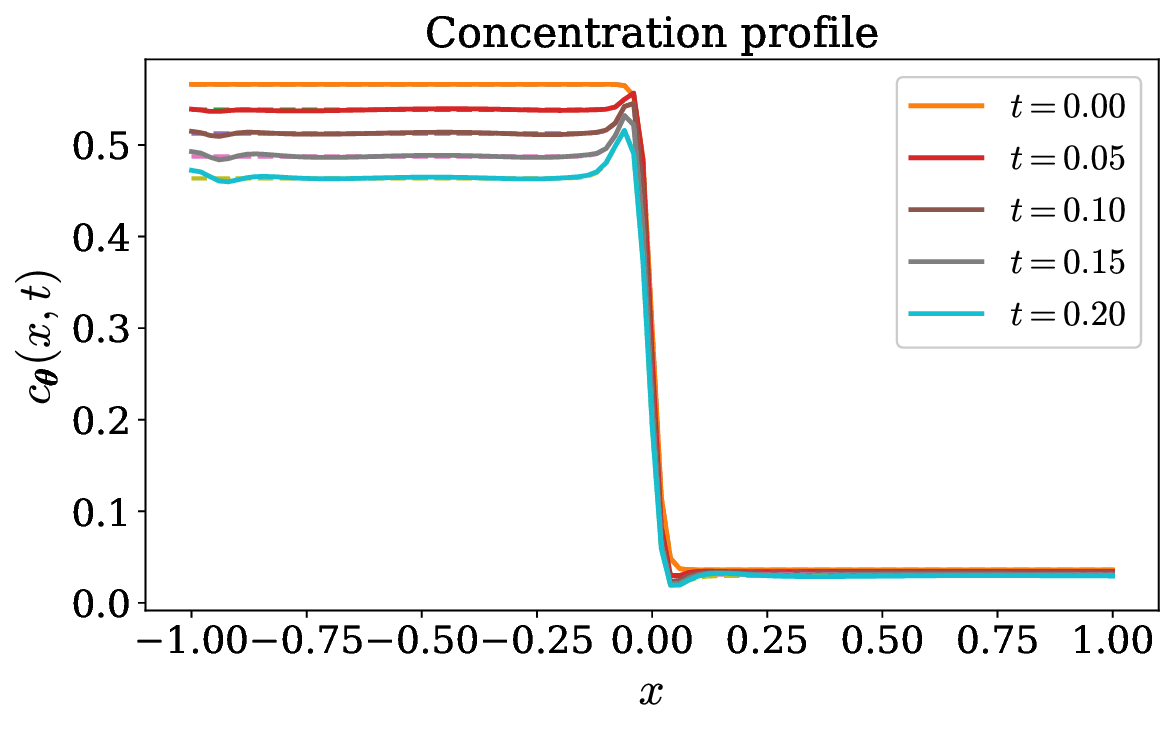}
    \caption{}
    \label{fig:ex2_diffusion_c_profiles}
  \end{subfigure}
  \caption{Profile comparison for Example~\ref{subsec:ex2_activation_diffusion}. Subfigures (a) and (b) show the phase-field profile $\phi_{\bth}(x,t)$ and concentration profile $c_{\bth}(x,t)$ in the activation-controlled regime. Subfigures (c), (d) show the corresponding profiles in the diffusion-controlled regime.}
  \label{fig:ex2_profiles}
\end{figure}

We now compare the performance of FM-tfPINN for the activation-controlled and diffusion-controlled regimes. 
The objective is to determine whether the same fractional-memory neural representation can resolve the coupled phase-field and concentration dynamics when the kinetic mechanism is changed through the phase-field mobility. Figure~\ref{fig:ex2_profiles} shows the predicted profiles for both regimes.
In each case, the learned phase field $\phi_{\bth}$ follows the reference diffuse interface, while the learned concentration $c_{\bth}$ reproduces the concentration redistribution across the interfacial region. 
The agreement is especially relevant near the moving front, where the phase transition and concentration gradient are strongly coupled.

The quantitative diagnostics are summarized in Table~\ref{tab:ex2_diagnostics}. In both regimes, the concentration field is recovered with relative $L^2$ accuracy of order $10^{-3}$, while the phase-field relative error remains at the order of $10^{-2}$. More importantly, the physical diagnostics $d(t)$, $c_{\Gamma}(t)$, and $c_{\max}(t)$ are reproduced with small relative errors in both regimes. This confirms that the learned solution retains the physically relevant corrosion observables, not only the pointwise field values.

\begin{table}[htbp]
\centering
\scriptsize
\caption{Quantitative diagnostics for Example~\ref{subsec:ex2_activation_diffusion}.}
\label{tab:ex2_diagnostics}
\renewcommand{\arraystretch}{1.22}
\begin{tabularx}{\textwidth}{@{}Xcc@{}}
\toprule
Diagnostic quantity & Activation & Diffusion \\
\midrule
Relative $L^2$ error of $\phi$, $E_2^\phi$ 
& $2.037352\mathrm{e}{-}02$ & $2.167774\mathrm{e}{-}02$ \\

Relative $L^2$ error of $c$, $E_2^c$ 
& $1.747944\mathrm{e}{-}03$ & $2.089513\mathrm{e}{-}03$ \\

Maximum error of $\phi$, $E_\infty^\phi$ 
& $2.748618\mathrm{e}{-}02$ & $1.963470\mathrm{e}{-}02$ \\

Maximum error of $c$, $E_\infty^c$ 
& $6.914768\mathrm{e}{-}03$ & $8.201043\mathrm{e}{-}03$ \\

Relative error of corrosion depth $d(t)$ 
& $3.412541\mathrm{e}{-}03$ & $3.467779\mathrm{e}{-}03$ \\

Relative error of interface concentration $c_\Gamma(t)$ 
& $4.317552\mathrm{e}{-}04$ & $7.307622\mathrm{e}{-}04$ \\

Relative error of concentration peak $c_{\max}(t)$ 
& $4.098336\mathrm{e}{-}04$ & $5.774832\mathrm{e}{-}04$ \\

Final predicted corrosion depth $d_{\bth}$ 
& $1.049582\mathrm{e}{-}02$ & $2.089889\mathrm{e}{-}02$ \\

Final reference corrosion depth $d_{\rm ref}$ 
& $1.050628\mathrm{e}{-}02$ & $2.081510\mathrm{e}{-}02$ \\

Final predicted interface concentration $c_{\Gamma,\bth}$ 
& $3.335405\mathrm{e}{-}01$ & $3.861993\mathrm{e}{-}01$ \\

Final reference interface concentration $c_{\Gamma,\rm ref}$ 
& $3.337007\mathrm{e}{-}01$ & $3.860398\mathrm{e}{-}01$ \\

Final predicted concentration peak $c_{\max,\bth}$ 
& $4.892262\mathrm{e}{-}01$ & $5.039443\mathrm{e}{-}01$ \\

Final reference concentration peak $c_{\max,\rm ref}$ 
& $4.895047\mathrm{e}{-}01$ & $5.039606\mathrm{e}{-}01$ \\

Maximum tempered concentration-mass drift 
& $2.715188\mathrm{e}{-}04$ & $7.776297\mathrm{e}{-}04$ \\
\bottomrule
\end{tabularx}
\renewcommand{\arraystretch}{1}
\end{table}

\begin{figure}[htbp]
  \centering
  \begin{subfigure}[b]{0.32\textwidth}
    \includegraphics[width=\textwidth]{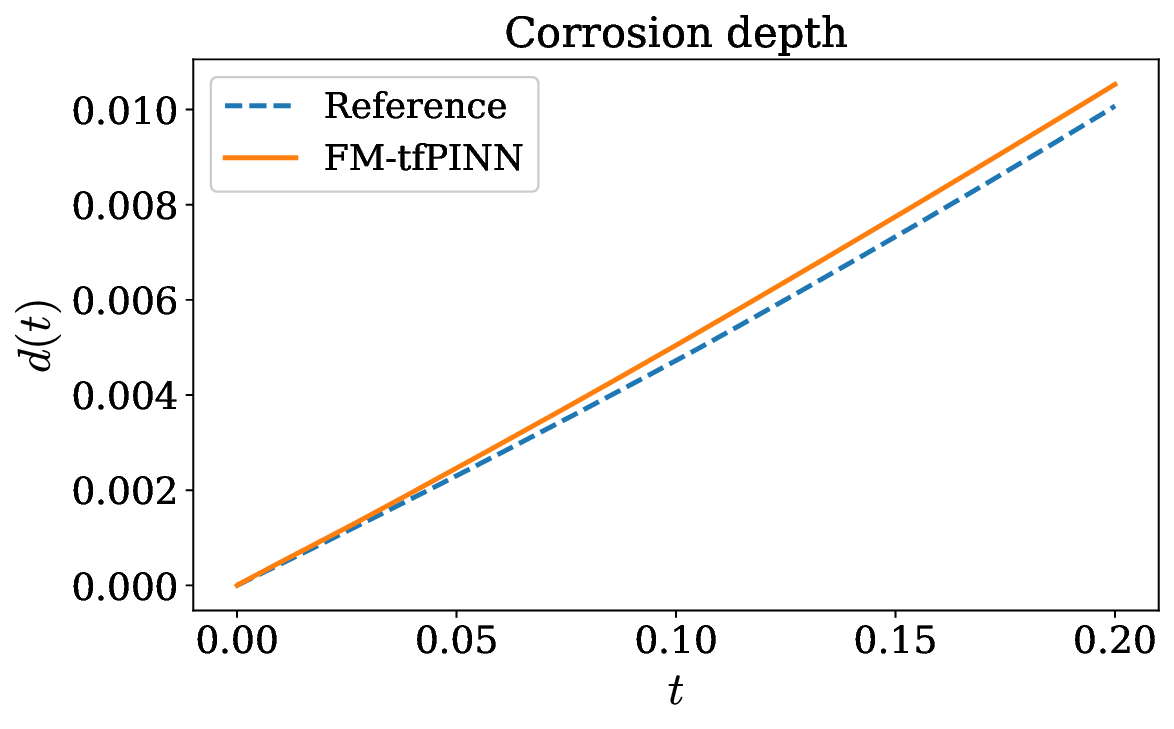}
    \caption{}
    \label{fig:ex2_activation_depth}
  \end{subfigure}
  \hfill
  \begin{subfigure}[b]{0.32\textwidth}
    \includegraphics[width=\textwidth]{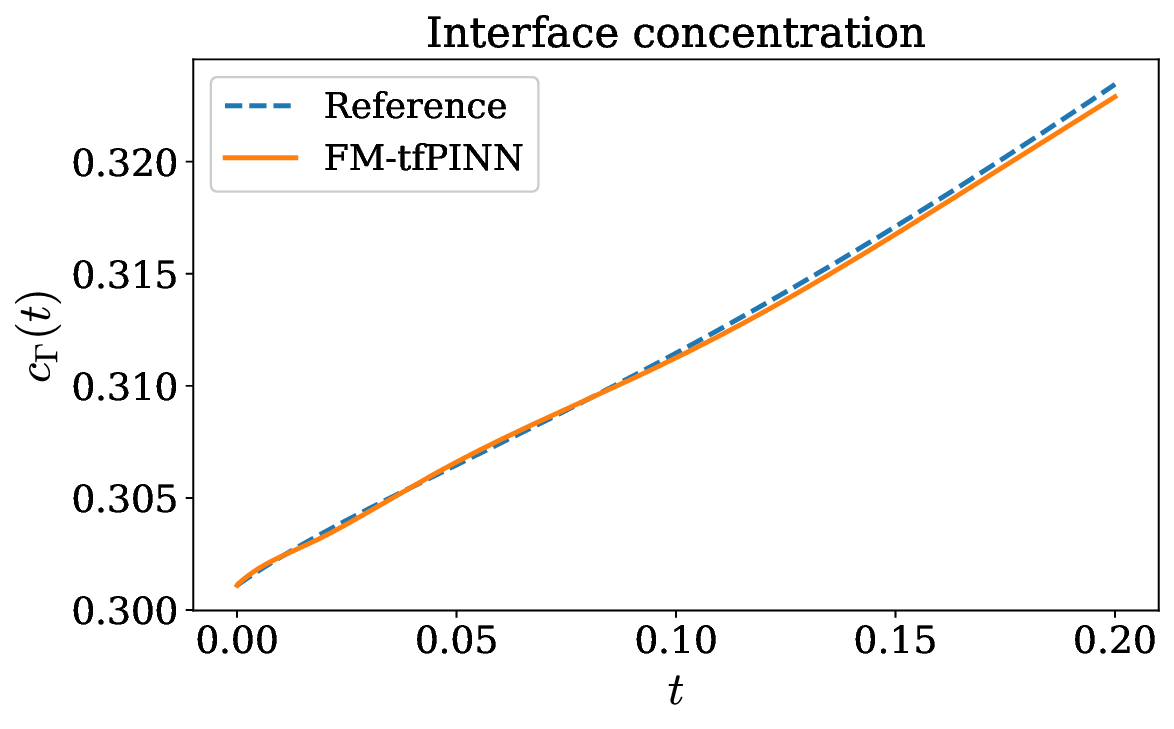}
    \caption{}
    \label{fig:ex2_activation_interface_c}
  \end{subfigure}
  \hfill
  \begin{subfigure}[b]{0.32\textwidth}
    \includegraphics[width=\textwidth]{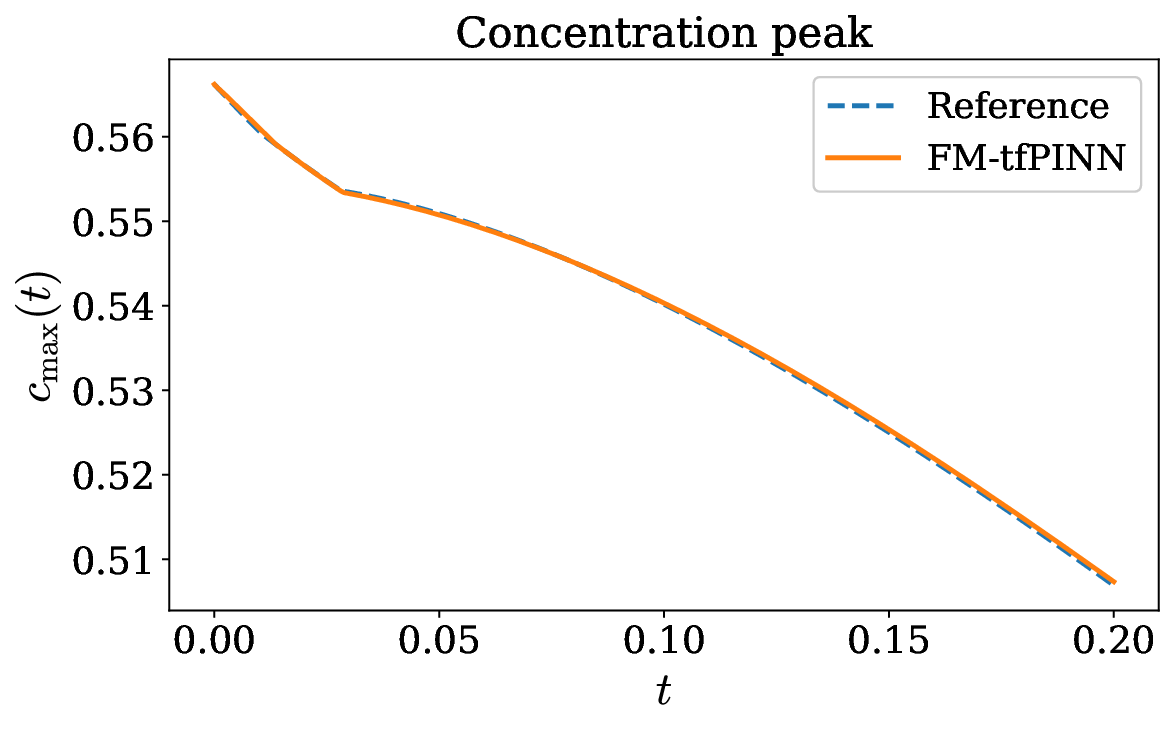}
    \caption{}
    \label{fig:ex2_activation_cmax}
  \end{subfigure}

  \begin{subfigure}[b]{0.32\textwidth}
    \includegraphics[width=\textwidth]{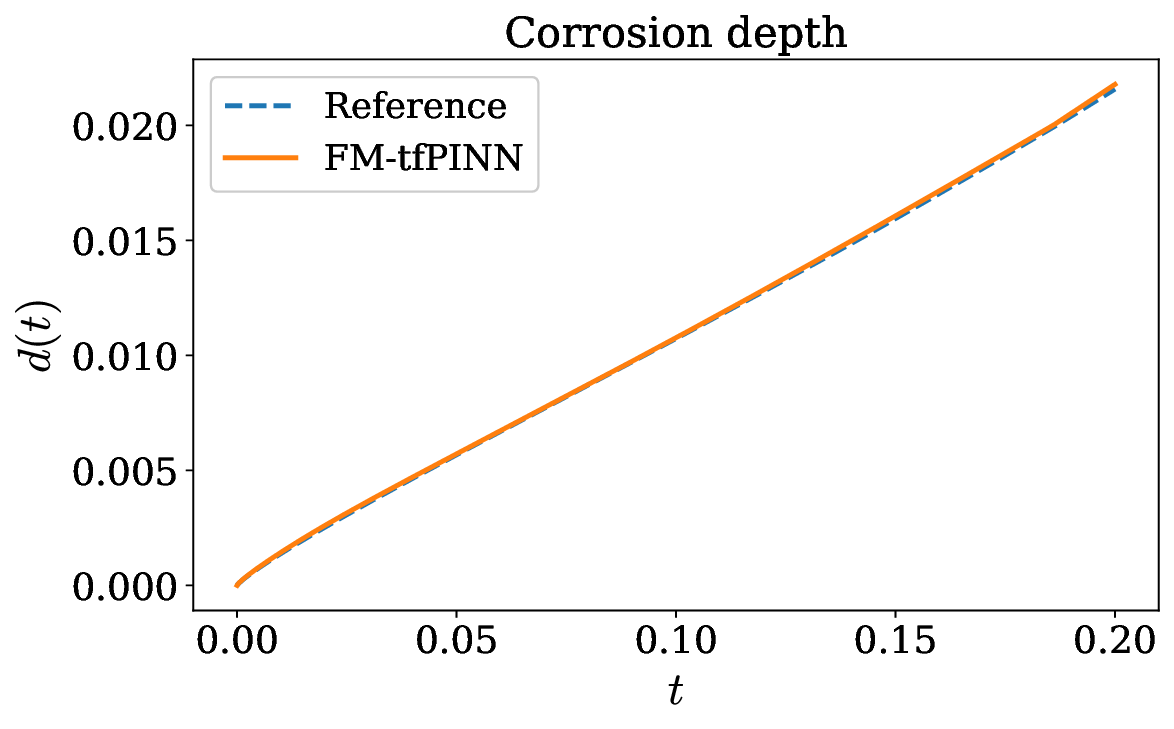}
    \caption{}
    \label{fig:ex2_diffusion_depth}
  \end{subfigure}
  \hfill
  \begin{subfigure}[b]{0.32\textwidth}
    \includegraphics[width=\textwidth]{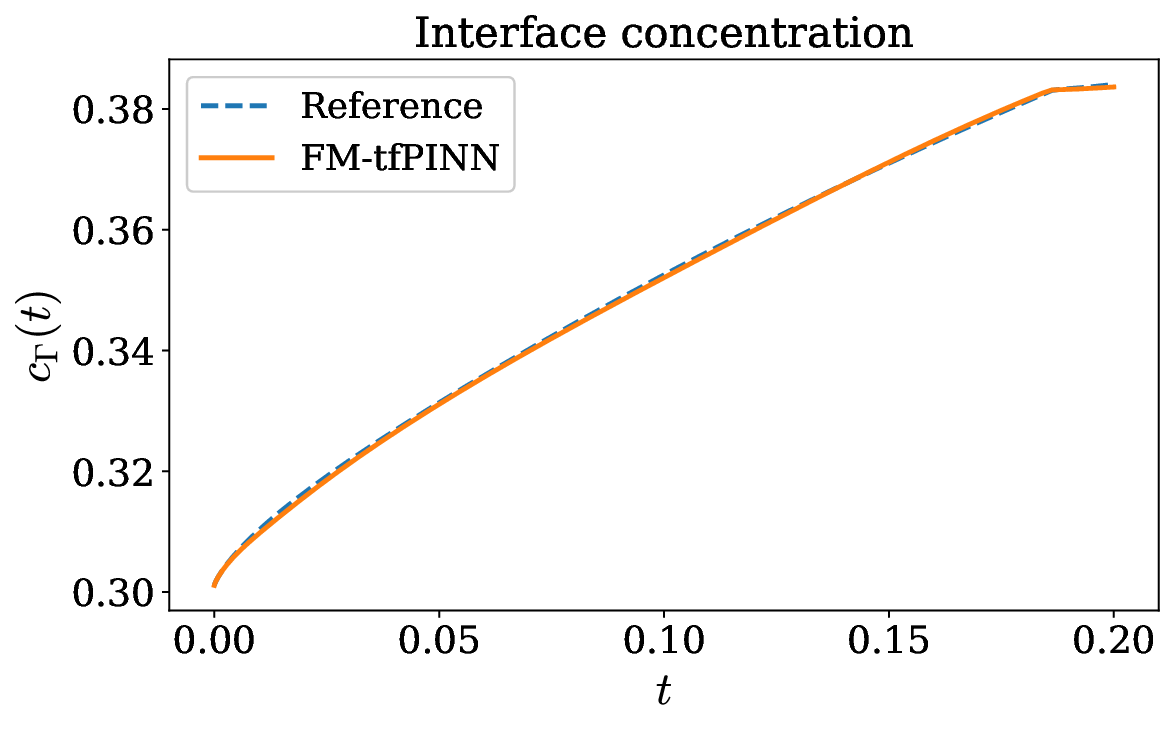}
    \caption{}
    \label{fig:ex2_diffusion_interface_c}
  \end{subfigure}
  \hfill
  \begin{subfigure}[b]{0.32\textwidth}
    \includegraphics[width=\textwidth]{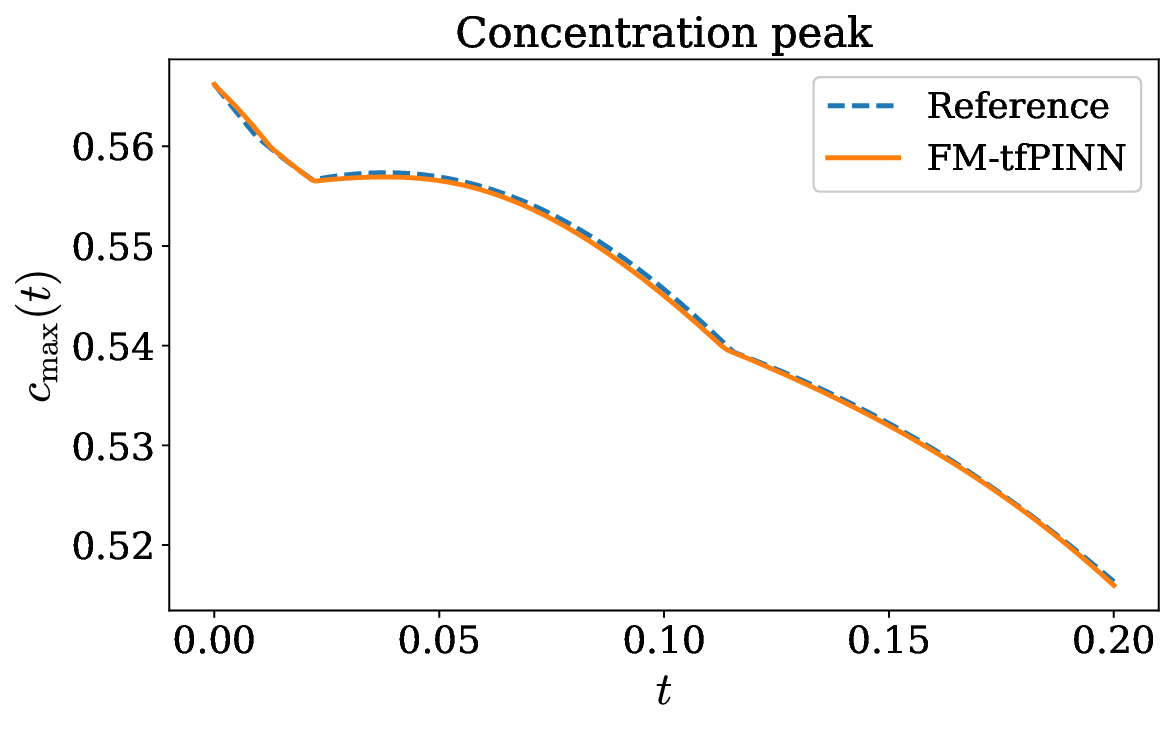}
    \caption{}
    \label{fig:ex2_diffusion_cmax}
  \end{subfigure}
  \caption{Physical diagnostics for Example~\ref{subsec:ex2_activation_diffusion}. Subfigures (a)--(c) show the corrosion depth $d(t)$, interface concentration $c_\Gamma(t)$, and concentration peak $c_{\max}(t)$ in the activation-controlled regime. Subfigures (d)--(f) show the corresponding diagnostics in the diffusion-controlled regime.}
  \label{fig:ex2_physical_diagnostics}
\end{figure}

The physical diagnostics are shown in Figure~\ref{fig:ex2_physical_diagnostics}. 
The corrosion depth $d_{\bth}(t)$ agrees closely with the reference curve in both regimes.
The diffusion-controlled case produces a larger corrosion depth than the activation-controlled case, which is consistent with the increased phase-field mobility and the faster interfacial advancement. 
The interface concentration $c_{\Gamma,\bth}(t)$ and concentration peak $c_{\max,\bth}(t)$ are also accurately reproduced. 
These quantities are sensitive to the coupling between interface motion and concentration transport, and therefore provide a stronger physical validation than field errors alone.

The comparison between the two regimes demonstrates that FM-tfPINN is not restricted to a single kinetic setting. With the same memory-generated representation and shifted residual construction, the method captures both the slower activation-controlled response and the faster diffusion-controlled front propagation. The small errors in the regime-sensitive diagnostics further show that the proposed framework can recover physically meaningful corrosion observables from the coupled tempered fractional dynamics.

\subsection{Semi-Circular Pitting Corrosion}\label{subsec:ex3_semicircular_pitting}

\subsubsection{Problem Statement}\label{subsubsec:ex3_problem_statement}
We next consider a two-dimensional semi-circular pitting corrosion problem. 
This example extends the preceding one-dimensional tests to a curved moving interface and is designed to assess whether FM-tfPINN can resolve localized pitting growth in a two-dimensional geometry.
The unknowns are again the phase-field variable $\phi(x,y,t)$ and the concentration field $c(x,y,t)$, where $\phi=1$ denotes the metal phase and $\phi=0$ denotes the electrolyte or pit region. 
For $0<\alpha<1$ and $\lambda>0$, the two-dimensional tTFCP system is given by
\begin{equation}\label{eq:ex3_pitting_model}
\begin{cases}
\partial_t^{\alpha,\lambda}\phi=2A_{\rm chem}L_{\rm mob}
   \bigl[c-h(\phi)(c_{\rm Se}-c_{\rm Le})-c_{\rm Le}\bigr]
    (c_{\rm Se}-c_{\rm Le})h'(\phi)\\
   \qquad
   -L_{\rm mob}w_{\phi}g'(\phi)+L_{\rm mob}\alpha_{\phi}\Delta\phi,\\
\partial_t^{\alpha,\lambda}c
=2A_{\rm chem}M_{\rm mob}\Delta c-2A_{\rm chem}M_{\rm mob}(c_{\rm Se}-c_{\rm Le})\Delta h(\phi),\\
\phi(x,y,0)=\phi_0(x,y),
\qquad
c(x,y,0)=c_0(x,y),
\qquad (x,y)\in\overline{\Omega},\\
\nabla \phi(x,y,t)\cdot\boldsymbol{n}=0,
\qquad
\nabla c(x,y,t)\cdot\boldsymbol{n}=0,
\qquad (x,y)\in\partial\Omega,\quad t\in[0,T],
\end{cases}
\end{equation}
where $\boldsymbol{n}$ denotes the outward unit normal on $\partial\Omega$. 
The interpolation function $h(\phi)$ and the double-well potential $g(\phi)$ are the same as those used in Example~\ref{subsec:ex1_pencil_electrode}. The initial pit is represented by a diffuse semi-circular interface attached to the lower boundary. 
Let $r_p(x,y)=\sqrt{(x-x_c)^2+(y-y_c)^2}$,
where $(x_c,y_c)$ is the center of the initial pit and $R_0$ is its initial radius. 
The initial phase field is prescribed by
\begin{equation*}
  \phi_0(x,y)
  = \frac{1}{2}\biggl[1+\tanh\Bigl(\frac{r_p(x,y)-R_0}{\sqrt{2}\,\ell}\Bigr)\biggr],
\qquad\ell=\sqrt{\frac{\alpha_{\phi}}{w_{\phi}}}.
\end{equation*}
Thus, the region inside the semi-circular pit corresponds to the electrolyte phase, while the exterior corresponds to the metal phase.
The initial concentration is chosen consistently with the phase interpolation,
\begin{equation*}
  c_0(x,y)=c_{\rm Le}+h(\phi_0(x,y))(c_{\rm Se}-c_{\rm Le}).
\end{equation*}
This experiment therefore complements the one-dimensional tests by evaluating FM-tfPINN in a genuinely two-dimensional interfacial setting, where the corrosion front is curved, spatially localized, and coupled to concentration transport in both coordinate directions. It provides a geometric validation of the proposed framework for resolving semi-circular pitting evolution governed by tempered fractional-memory dynamics.

\subsubsection{Implementation Details}\label{subsubsec:ex3_implementation_details}
The computation is performed according to the FM-tfPINN training procedure summarized in Algorithm~\ref{alg:fm_tfpinn}. 
In contrast to Example~\ref{subsec:ex1_pencil_electrode}, the present two-dimensional pitting problem involves homogeneous no-flux conditions for both coupled variables.
Hence, no endpoint-vanishing Dirichlet mask is introduced for the phase-field component. 
The memory-generated approximation is therefore written as
\begin{equation}\label{eq:ex3_fm_representation}
   \phi_{\bth}(x,y,t)
   =\phi_0(x,y)+\mathcal{I}_t^{\alpha,\lambda}[z_{\phi,\bth}](x,y,t),
\quad
c_{\bth}(x,y,t)=c_0(x,y)+\mathcal{I}_t^{\alpha,\lambda}[z_{c,\bth}](x,y,t).
\end{equation}
The initial state is embedded in the approximation by the fractional-memory representation, since the tempered fractional integral vanishes at the initial time. 
The homogeneous no-flux boundary conditions are enforced weakly through the boundary-loss terms, and the data-mismatch contribution is omitted because this test is performed in the forward-prediction setting. 
A reference solution is computed using the SOE-accelerated shifted finite-difference discretization employed in the preceding tests, extended here to the two-dimensional system \eqref{eq:ex3_pitting_model}.
The FM-tfPINN residual is evaluated on the shifted temporal levels $t_{n+\sigma}$, with $\sigma=1-\alpha/2$, and on spatial collocation points distributed in $\Omega$. 

To enhance the resolution of the localized pitting front, the residual set is enriched by points sampled from the initial diffuse-interface band and from the predicted interfacial region
$\bigl\{(x,y)\in\Omega\colon\phi_{\min}<\phi_{\bth}(x,y,t)<\phi_{\max}\bigr\}$.
Additional residual-adaptive points are selected from regions where the coupled residuals are relatively large.
This sampling strategy is well suited to the present 2D 
pitting problem, where the dominant dynamics are concentrated around a curved moving front, while a large part of the domain remains close to a slowly varying bulk state. 
Besides the field errors $E_2^\phi$, $E_2^c$, $E_\infty^\phi$, and $E_\infty^c$, we use two geometry-based diagnostics to quantify the learned pitting evolution. 
The pit area is defined by
\begin{equation}\label{eq:ex3_pit_area}
   A_{\rm pit}(t)=\Bigl|\bigl\{(x,y)\in\Omega\colon\phi(x,y,t)<0.5\bigr\}\Bigr|,
\end{equation}
and the corresponding equivalent semi-circular radius is computed as
\begin{equation*}
    R_{\rm eq}(t)=\sqrt{\frac{2A_{\rm pit}(t)}{\pi}}.
\end{equation*}
These quantities provide physically interpretable measures of pit growth and allow the accuracy of the learned two-dimensional interface motion to be assessed beyond pointwise field errors.
The implementation settings are reported in Table~\ref{tab:ex3_implementation_details}.

\begin{table}[htbp]
\centering
\scriptsize
\caption{Implementation details for Example~\ref{subsec:ex3_semicircular_pitting}.}
\label{tab:ex3_implementation_details}
\renewcommand{\arraystretch}{1.25}
\begin{tabular}{ll}
\toprule
Item & Setting \\
\midrule
Configuration & FM-tfPINN forward prediction for two-dimensional pitting \\
Computational domain & $\Omega=(-1,1)\times(0,1)$ \\
Initial geometry & Semi-circular diffuse pit attached to the lower boundary \\
Fractional order and tempering & $\alpha=0.8,\quad \lambda=1$ \\
Shifted level & $\sigma=1-\alpha/2=0.6$ \\
FM-tfPINN temporal grid & $N_t=20,\quad r=(3-\alpha)/\alpha$ \\
Boundary condition & Homogeneous no-flux conditions for $\phi$ and $c$ \\
Neural representation & $(x,y,t)\mapsto(z_{\phi,\bth},z_{c,\bth})$ through \eqref{eq:ex3_fm_representation} \\
Collocation strategy & Bulk, initial-interface, self-interface, and residual-adaptive points \\
Self-interface band & $\phi_{\min}<\phi_{\bth}<\phi_{\max}$ \\
Physical diagnostics & Pit area $A_{\rm pit}(t)$ and equivalent radius $R_{\rm eq}(t)$ \\
Observation data & Not used, $\omega_d=0$ \\
\bottomrule
\end{tabular}
\renewcommand{\arraystretch}{1}
\end{table}

\subsubsection{Results and Discussion}\label{subsubsec:ex3_results_discussion}
\begin{figure}[htbp]
  \centering
  \begin{subfigure}[b]{0.96\textwidth}
    \includegraphics[width=\textwidth]{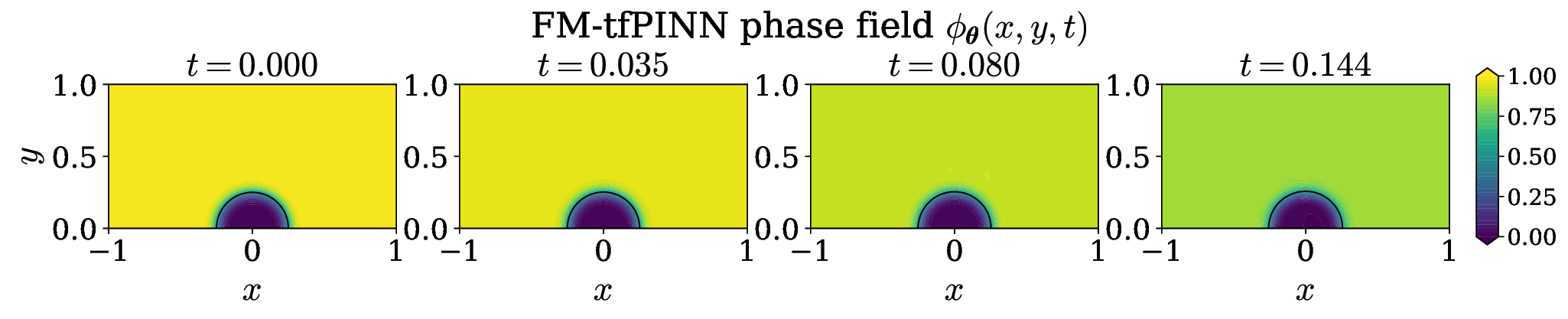}
    \caption{}
    \label{fig:ex3_phi_contours}
  \end{subfigure}
  \begin{subfigure}[b]{0.96\textwidth}
    \includegraphics[width=\textwidth]{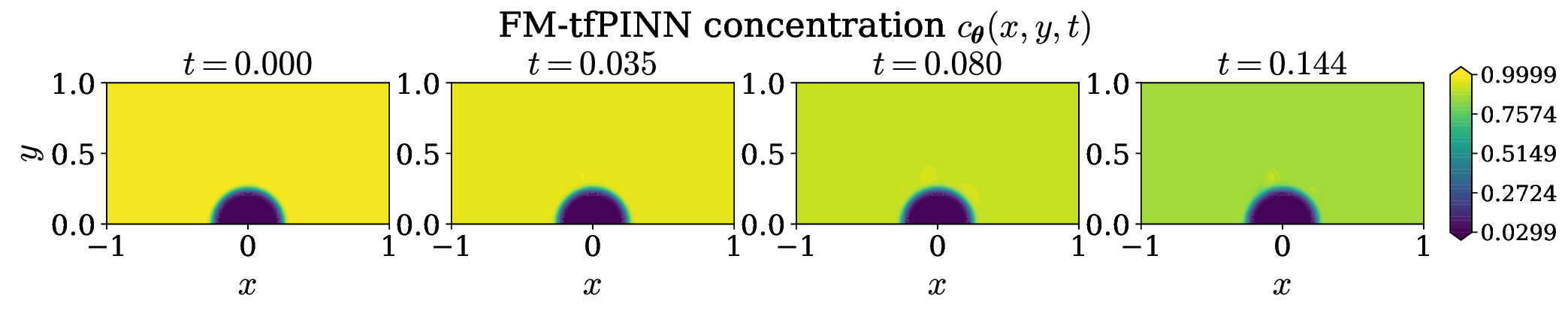}
    \caption{}
    \label{fig:ex3_c_contours}
  \end{subfigure}
  \caption{FM-tfPINN prediction for Example~\ref{subsec:ex3_semicircular_pitting}.
  Subfigure (a) shows the phase-field distribution $\phi_{\bth}(x,y,t)$. 
  Subfigure (b) shows the concentration distribution $c_{\bth}(x,y,t)$.}
  \label{fig:ex3_contours}
\end{figure}

We now evaluate the proposed FM-tfPINN framework for the 2D 
semi-circular pitting corrosion problem. This example is intended to test whether the memory-generated approximation can resolve a localized curved interface and recover geometry-dependent corrosion observables under the tempered fractional dynamics.
Figure~\ref{fig:ex3_contours} shows the learned phase-field and concentration distributions at representative shifted time levels. 
The predicted phase field $\phi_{\bth}(x,y,t)$ preserves the semi-circular pit morphology and captures the evolution of the curved diffuse interface. 
The concentration field $c_{\bth}(x,y,t)$ evolves consistently with the pitting region, indicating that the coupled interfacial and concentration dynamics are learned simultaneously in the 2D 
setting.

\begin{figure}[htbp]
  \centering
  \begin{subfigure}[b]{0.96\textwidth}
    \includegraphics[width=\textwidth]{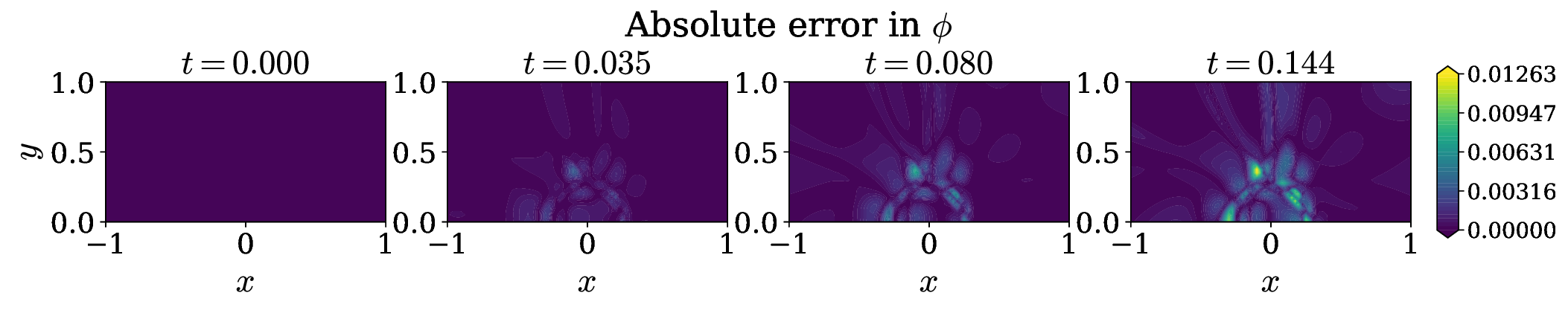}
    \caption{}
    \label{fig:ex3_phi_error_contours}
  \end{subfigure}
  \begin{subfigure}[b]{0.96\textwidth}
    \includegraphics[width=\textwidth]{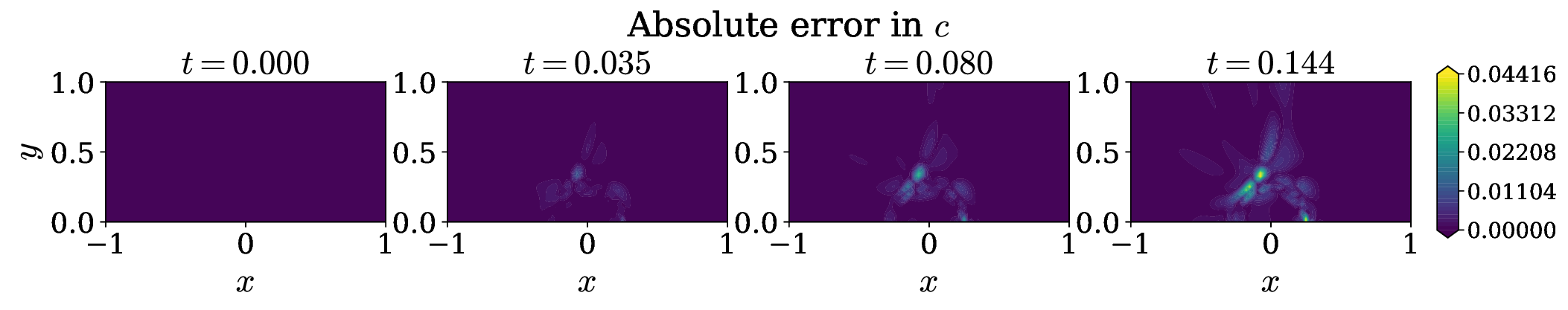}
    \caption{}
    \label{fig:ex3_c_error_contours}
  \end{subfigure}
  \caption{Absolute-error contours for Example~\ref{subsec:ex3_semicircular_pitting}. Subfigure (a) shows $|\phi_{\bth}-\phi_{\rm ref}|$. 
  Subfigure (b) shows $|c_{\bth}-c_{\rm ref}|$.}
  \label{fig:ex3_error_contours}
\end{figure}

The corresponding absolute-error distributions are displayed in Figure~\ref{fig:ex3_error_contours}. The phase-field error remains small and is concentrated around the evolving semi-circular interface, which is the most sensitive region of the computation. The concentration error is also localized primarily near the pit boundary and the associated concentration-gradient region. 
This behavior is expected for a moving-interface phase-field problem, since small geometric shifts in the diffuse front produce localized pointwise deviations even when the global field errors remain small.
The error contours therefore provide additional evidence that FM-tfPINN resolves the two-dimensional pitting morphology without introducing spurious global errors across the bulk region.

The quantitative diagnostics are summarized in Table~\ref{tab:ex3_diagnostics}.
The relative field errors remain small for both components, with the phase-field error below $10^{-3}$ and the concentration error of order $10^{-3}$. 
More importantly, the final pit area and equivalent radius computed from the learned phase field agree with the reference values on the shifted evaluation grid. 
This confirms that FM-tfPINN recovers not only the field variables, but also the geometric quantities that characterize the growth of the pitting region.

\begin{table}[htbp]
\centering
\scriptsize
\caption{Quantitative diagnostics for Example~\ref{subsec:ex3_semicircular_pitting}.}
\label{tab:ex3_diagnostics}
\renewcommand{\arraystretch}{1.25}
\begin{tabular}{lc}
\toprule
Diagnostic quantity & Value \\
\midrule
Relative $L^2$ error of the phase field, $E_2^\phi$ & $7.043799\mathrm{e}{-}04$ \\
Relative $L^2$ error of the concentration field, $E_2^c$ & $1.508610\mathrm{e}{-}03$ \\
Maximum error of the phase field, $E_\infty^\phi$ & $1.262852\mathrm{e}{-}02$ \\
Maximum error of the concentration field, $E_\infty^c$ & $1.415560\mathrm{e}{-}02$ \\
Predicted final pit area, $A_{{\rm pit},\bth}(t_{N_t-1+\sigma})$ & $1.106250\mathrm{e}{-}01$ \\
Reference final pit area, $A_{{\rm pit},{\rm ref}}(t_{N_t-1+\sigma})$ & $1.106250\mathrm{e}{-}01$ \\
Predicted final equivalent radius, $R_{{\rm eq},\bth}(t_{N_t-1+\sigma})$ & $2.653791\mathrm{e}{-}01$ \\
Reference final equivalent radius, $R_{{\rm eq},{\rm ref}}(t_{N_t-1+\sigma})$ & $2.653791\mathrm{e}{-}01$ \\
\bottomrule
\end{tabular}
\renewcommand{\arraystretch}{1}
\end{table}

\begin{figure}[htbp]
  \centering
  \begin{subfigure}[b]{0.48\textwidth}
    \includegraphics[width=\textwidth]{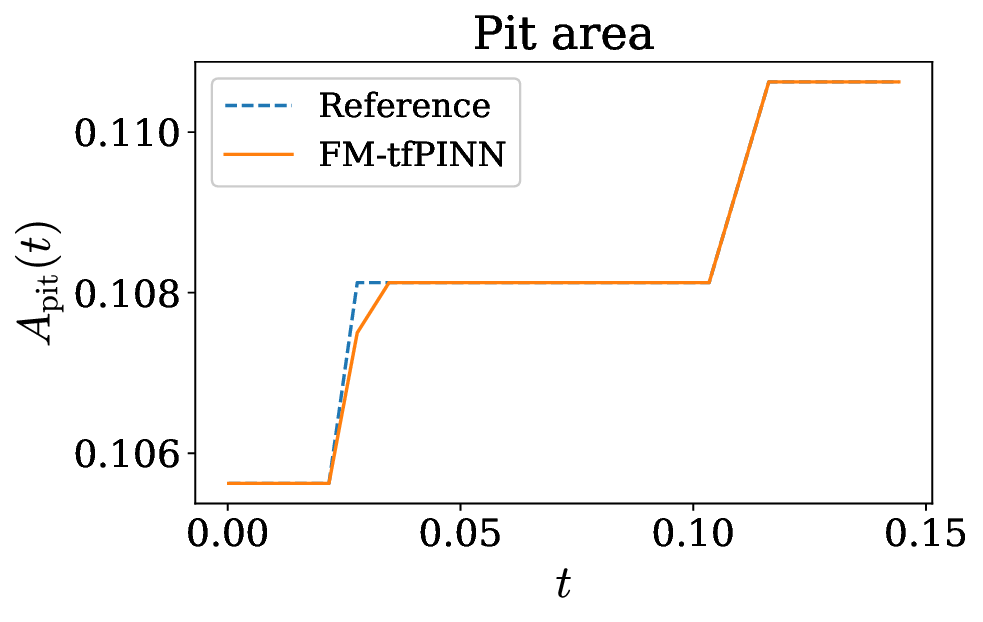}
    \caption{}
    \label{fig:ex3_pit_area}
  \end{subfigure}
  \hfill
  \begin{subfigure}[b]{0.48\textwidth}
    \includegraphics[width=\textwidth]{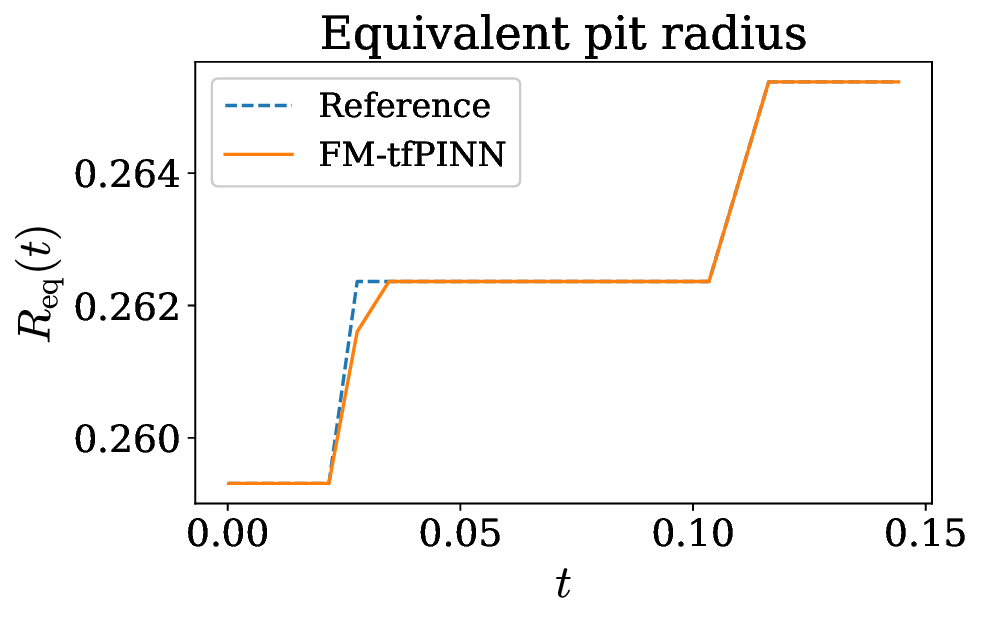}
    \caption{}
    \label{fig:ex3_equivalent_radius}
  \end{subfigure}
  \caption{Geometry-based pitting diagnostics for Example~\ref{subsec:ex3_semicircular_pitting}. Subfigure (a) compares the pit area $A_{\rm pit}(t)$. Subfigure (b) compares the equivalent pit radius $R_{\rm eq}(t)$.}
  \label{fig:ex3_pitting_diagnostics}
\end{figure}

Figure~\ref{fig:ex3_pitting_diagnostics} further compares the physical pitting diagnostics obtained from FM-tfPINN with the reference solution. 
The predicted pit area $A_{\rm pit}(t)$ follows the reference curve throughout the simulation interval. 
The equivalent radius $R_{\rm eq}(t)$, which converts the learned pit area into a geometry-based radius measure, also agrees with the reference result. 
Since both diagnostics are extracted from the learned diffuse interface and are not directly imposed as training constraints, their agreement provides a strong validation of the proposed interface-aware and residual-adaptive FM-tfPINN formulation.

Overall, this example demonstrates that FM-tfPINN can accurately resolve two-dimensional pitting corrosion with a curved moving interface. 
The agreement of the learned phase-field and concentration distributions, the localized structure of the absolute errors, and the recovery of the pit area and equivalent radius show that the proposed framework is effective not only for one-dimensional corrosion-front propagation but also for localized two-dimensional interface evolution governed by tempered fractional memory.

\subsection{Inverse Identification from Physical Corrosion Observations}\label{subsec:ex4_inverse_identification}

\subsubsection{Problem Statement}\label{subsubsec:ex4_problem_statement}
We finally consider an inverse problem for identifying hidden kinetic parameters in the tTFCP system from sparse physical corrosion observations. 
The computational domain is $\Omega=(-H_s,H_l)$, and the coupled variables are the phase-field variable $\phi(x,t)$ and the normalized concentration $c(x,t)$. 
In contrast to Example~\ref{subsec:ex1_pencil_electrode}, the present inverse setting is formulated with homogeneous no-flux boundary conditions for both variables. 
The governing problem is
\begin{equation}\label{eq:ex4_inverse_model}
\begin{cases}
\partial_t^{\alpha,\lambda}\phi
=2A_{\rm chem}L_{\rm mob} \bigl[c-h(\phi)(c_{\rm Se}-c_{\rm Le})-c_{\rm Le}\bigr]
(c_{\rm Se}-c_{\rm Le})h'(\phi)
\\ \qquad\qquad  -L_{\rm mob}w_{\phi}g'(\phi)+L_{\rm mob}\alpha_{\phi}\phi_{xx},\\
\partial_t^{\alpha,\lambda}c 
=2A_{\rm chem}M_{\rm mob}c_{xx}-2A_{\rm chem}M_{\rm mob}(c_{\rm Se}-c_{\rm Le})(h(\phi))_{xx},\\
\phi(x,0)=\phi_0(x),\qquad c(x,0)=c_0(x),\qquad x\in[-H_s,H_l],\\
\phi_x(-H_s,t) =\phi_x(H_l,t)=0,
\qquad c_x(-H_s,t)=c_x(H_l,t)=0,
\qquad t\in[0,T].
\end{cases}
\end{equation}
The interpolation function $h(\phi)$, the double-well potential $g(\phi)$, and the interface thickness $\ell$ are the same as in the preceding forward examples. The initial phase field is chosen as
\begin{equation*} 
\phi_0(x)=\frac{1}{2}\biggl[1-\tanh\Bigl(\frac{x}{\sqrt{2}\,\ell}\Bigr)\biggr],
\qquad\ell=\sqrt{\frac{\alpha_{\phi}}{w_{\phi}}}.
\end{equation*}
To generate an inverse setting based on concentration-deficient corrosion data, the initial concentration is prescribed in the undersaturated form
\begin{equation*}
c_0(x) = c_{\rm Le}+ \eta_0 h(\phi_0(x))(c_{\rm Se}-c_{\rm Le}),
\qquad \eta_0=0.55 .
\end{equation*}
The inverse task is to recover the mobility vector
$\bka=(L_{\rm mob},M_{\rm mob})$
from a sparse set of physical observations. The true parameters used to generate the reference data are denoted by
$\bka^{\dagger}=(L_{\rm mob}^{\dagger},M_{\rm mob}^{\dagger})$.
The memory parameters $\alpha$ and $\lambda$ are kept fixed in this inverse test. This choice isolates the recovery of the kinetic and transport mobilities, which directly influence the corrosion-front motion and concentration redistribution. The observation data are not taken from dense solution snapshots. Instead, they consist of sparse physical corrosion measurements at selected observation times
$\mathcal{T}_{\rm obs}=\{t_{m_j}\}_{j=1}^{N_{\rm obs}}\subset(0,T]$.
The observation set is
\begin{equation}\label{eq:ex4_observation_set}
\mathcal{D}_{\rm obs}
=\Bigl\{\bigl(t_{m_j},d_{\rm obs}(t_{m_j}),
c_{\Gamma,{\rm obs}}(t_{m_j}),
c_{\max,{\rm obs}}(t_{m_j})\bigr)\Bigr\}_{j=1}^{N_{\rm obs}},
\end{equation}
where $d(t)$ is the corrosion depth extracted from the diffuse interface, $c_{\Gamma}(t)$ is the concentration near the interface and $c_{\max}(t)$ is the concentration peak. This observation design leads to a physically interpretable inverse problem, in which the unknown mobility parameters are identified from measurable corrosion diagnostics rather than from dense full-field solution data.

\subsubsection{Implementation Details}\label{subsubsec:ex4_implementation_details}
The computation follows Algorithm~\ref{alg:fm_tfpinn} with the inverse objective $\mathcal{L}_I$. Since the present problem uses homogeneous no-flux conditions for both $\phi$ and $c$, the endpoint-vanishing phase-field mask used in Example~\ref{subsec:ex1_pencil_electrode} is not introduced. 
The memory-generated approximation is therefore written as
\begin{equation}\label{eq:ex4_fm_representation}
\phi_{\bth}(x,t)
=\phi_0(x)+\mathcal{I}_t^{\alpha,\lambda}[z_{\phi,\bth}](x,t),
\qquad c_{\bth}(x,t)
=c_0(x)+\mathcal{I}_t^{\alpha,\lambda}[z_{c,\bth}](x,t).
\end{equation}
The initial state is incorporated through the vanishing of the tempered fractional integral at $t=0$. The no-flux boundary conditions are imposed weakly through the boundary-loss terms for both components. The unknown mobilities are treated as positive trainable parameters. To preserve positivity during optimization, we use the logarithmic parametrization
\begin{equation*} 
  L_{\rm mob}=\mathrm{e}^{\eta_L},\qquad M_{\rm mob}=\mathrm{e}^{\eta_M},
\end{equation*}
where $\eta_L$ and $\eta_M$ are optimized together with the neural-network parameters $\bth$. The full trainable set is therefore
$(\bth,\eta_L,\eta_M)$.
A reference solution is generated using the SOE-accelerated shifted finite-difference solver applied to \eqref{eq:ex4_inverse_model} with the true parameters $\bka^{\dagger}$. The sparse observations in $\mathcal{D}_{\rm obs}$ are extracted from this reference solution. 
During inverse training, only the physical observations $d_{\rm obs}$, $c_{\Gamma,{\rm obs}}$, and $c_{\max,{\rm obs}}$ are used in the data-mismatch loss.
The full reference fields are retained only for post-training validation. The data loss is defined by
\begin{equation}\label{eq:ex4_data_loss}
\mathcal{L}_{\rm data}
=\omega_d\mathcal{L}_d+\omega_{\Gamma}\mathcal{L}_{\Gamma}+\omega_{\max}\mathcal{L}_{\max},
\end{equation}
where
\begin{equation}\label{eq:ex4_observation_losses}
   \mathcal{L}_d=\frac{1}{N_{\rm obs}}
    \sum_{j=1}^{N_{\rm obs}}\bigl|d_{\bth}(t_{m_j})-d_{\rm obs}(t_{m_j})\bigr|^2,
\end{equation}
\begin{equation}\label{eq:ex4_interface_loss}
   \mathcal{L}_{\Gamma}=\frac{1}{N_{\rm obs}}
  \sum_{j=1}^{N_{\rm obs}}
  \bigl|c_{\Gamma,\bth}(t_{m_j})-c_{\Gamma,{\rm obs}}(t_{m_j})\bigr|^2,
\end{equation}
and
\begin{equation}\label{eq:ex4_cmax_loss}
  \mathcal{L}_{\max}=\frac{1}{N_{\rm obs}}
\sum_{j=1}^{N_{\rm obs}}
\bigl|c_{\max,\bth}(t_{m_j})-c_{\max,{\rm obs}}(t_{m_j})\bigr|^2 .
\end{equation}
The inverse objective is then $\mathcal{L}_I=\mathcal{L}_{F}+\mathcal{L}_{\rm data}$.
In the implementation, differentiable approximations of $d_{\bth}(t)$, $c_{\Gamma,\bth}(t)$, and $c_{\max,\bth}(t)$ are used during training so that the physical observation loss can be optimized by automatic differentiation. The main settings for this inverse experiment are summarized in Table~\ref{tab:ex4_implementation_details}.

\begin{table}[htbp]
\centering
\scriptsize
\caption{Implementation details for Example~\ref{subsec:ex4_inverse_identification}.}
\label{tab:ex4_implementation_details}
\renewcommand{\arraystretch}{1.25}
\begin{tabular}{ll}
\toprule
\textbf{Item} & \textbf{Setting} \\
\midrule
Configuration & FM-tfPINN inverse identification \\
Computational domain & $\Omega=(-1,1)$, with $H_s=H_l=1$ \\
Boundary condition & Homogeneous no-flux conditions for $\phi$ and $c$ \\
Fixed memory parameters & $\alpha=0.8,\quad \lambda=1$ \\
Unknown parameter vector & $\bka=(L_{\rm mob},M_{\rm mob})$ \\
True parameters & $L_{\rm mob}^{\dagger}=10^{-3},\quad M_{\rm mob}^{\dagger}=10^{-3}$ \\
Trainable parametrization & $L_{\rm mob}=\mathrm{e}^{\eta_L},\quad M_{\rm mob}=\mathrm{e}^{\eta_M}$ \\
Initial concentration mode & Undersaturated, with $\eta_0=0.55$ \\
Observation data & $d_{\rm obs}(t_m),\ c_{\Gamma,{\rm obs}}(t_m),\ c_{\max,{\rm obs}}(t_m)$ \\
Number of observation times & $N_{\rm obs}=10$ \\
Reference solver & SOE-accelerated shifted finite difference \\
FM-tfPINN representation & $(x,t)\mapsto(z_{\phi,\bth},z_{c,\bth})$ through \eqref{eq:ex4_fm_representation} \\
Training objective & $\mathcal{L}_I=\mathcal{L}_F+\mathcal{L}_{\rm data}$ \\
\bottomrule
\end{tabular}
\renewcommand{\arraystretch}{1}
\end{table}

\subsubsection{Results and Discussion}\label{subsubsec:ex4_results_discussion}
\begin{figure}[htbp]
  \centering
  \includegraphics[width=0.64\textwidth]{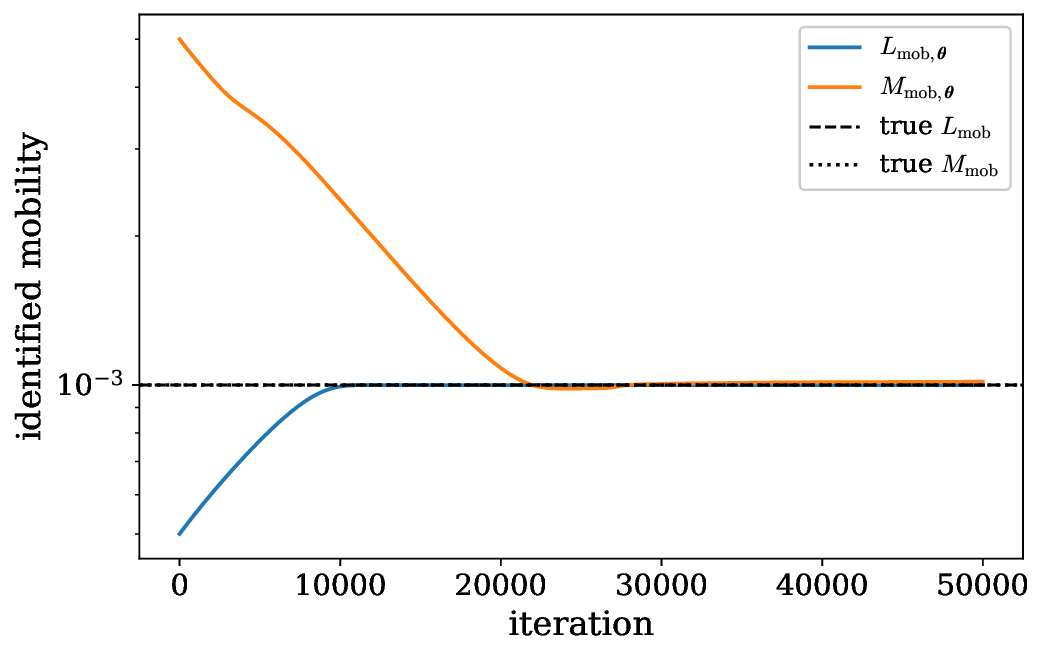}
  \caption{Parameter-identification history for Example~\ref{subsec:ex4_inverse_identification}. The curves show the learned mobilities $L_{{\rm mob},\bth}$ and $M_{{\rm mob},\bth}$, while the horizontal reference lines indicate $L_{\rm mob}^{\dagger}$ and $M_{\rm mob}^{\dagger}$.}
  \label{fig:ex4_parameter_convergence}
\end{figure}

We now assess the inverse capability of FM-tfPINN for identifying the mobility parameters in the tTFCP system from sparse physical corrosion observations. 
The inverse training uses only scalar corrosion observations at a small number of time instances, whereas the full reference fields are used only for post-training validation. Figure~\ref{fig:ex4_parameter_convergence} shows the convergence history of the identified mobility parameters. Both $L_{{\rm mob},\bth}$ and $M_{{\rm mob},\bth}$ converge toward their reference values during training. The phase-field mobility is recovered with very high accuracy, while the concentration mobility is also identified with a small relative deviation.

The quantitative inverse diagnostics are reported in Table~\ref{tab:ex4_inverse_diagnostics}.
The identified phase-field mobility is nearly identical to the true value, and the concentration mobility is recovered with a relative error of order $10^{-2}$. 
The field-reconstruction errors are also small, particularly for the phase-field component. 
The physical diagnostic errors for the corrosion depth, interface concentration, and concentration peak remain at a moderate level, which is acceptable for an inverse problem trained from sparse scalar observations rather than from full-field data. 

\begin{table}[htbp]
\centering
\scriptsize
\caption{Inverse-identification diagnostics for Example~\ref{subsec:ex4_inverse_identification}.}
\label{tab:ex4_inverse_diagnostics}
\renewcommand{\arraystretch}{1.25}
\begin{tabular}{lc}
\toprule
Diagnostic quantity & Value \\
\midrule
True phase-field mobility, $L_{\rm mob}^{\dagger}$ & $1.000000\mathrm{e}{-}03$ \\
Identified phase-field mobility, $L_{{\rm mob},\bth}$ & $9.999741\mathrm{e}{-}04$ \\
Relative error in $L_{\rm mob}$ & $2.594749\mathrm{e}{-}05$ \\
True concentration mobility, $M_{\rm mob}^{\dagger}$ & $1.000000\mathrm{e}{-}03$ \\
Identified concentration mobility, $M_{{\rm mob},\bth}$ & $1.015251\mathrm{e}{-}03$ \\
Relative error in $M_{\rm mob}$ & $1.525107\mathrm{e}{-}02$ \\
Relative $L^2$ error of the phase field, $E_2^\phi$ & $1.066528\mathrm{e}{-}03$ \\
Relative $L^2$ error of the concentration field, $E_2^c$ & $9.578217\mathrm{e}{-}03$ \\
Maximum error of the phase field, $E_\infty^\phi$ & $4.935226\mathrm{e}{-}03$ \\
Maximum error of the concentration field, $E_\infty^c$ & $2.218152\mathrm{e}{-}02$ \\
Relative error of the corrosion depth, $d(t)$ & $1.660202\mathrm{e}{-}02$ \\
Relative error of the interface concentration, $c_\Gamma(t)$ & $7.476640\mathrm{e}{-}03$ \\
Relative error of the concentration peak, $c_{\max}(t)$ & $2.468731\mathrm{e}{-}02$ \\
\bottomrule
\end{tabular}
\renewcommand{\arraystretch}{1}
\end{table}

\begin{figure}[htbp]
  \centering
  \includegraphics[width=0.96\textwidth]{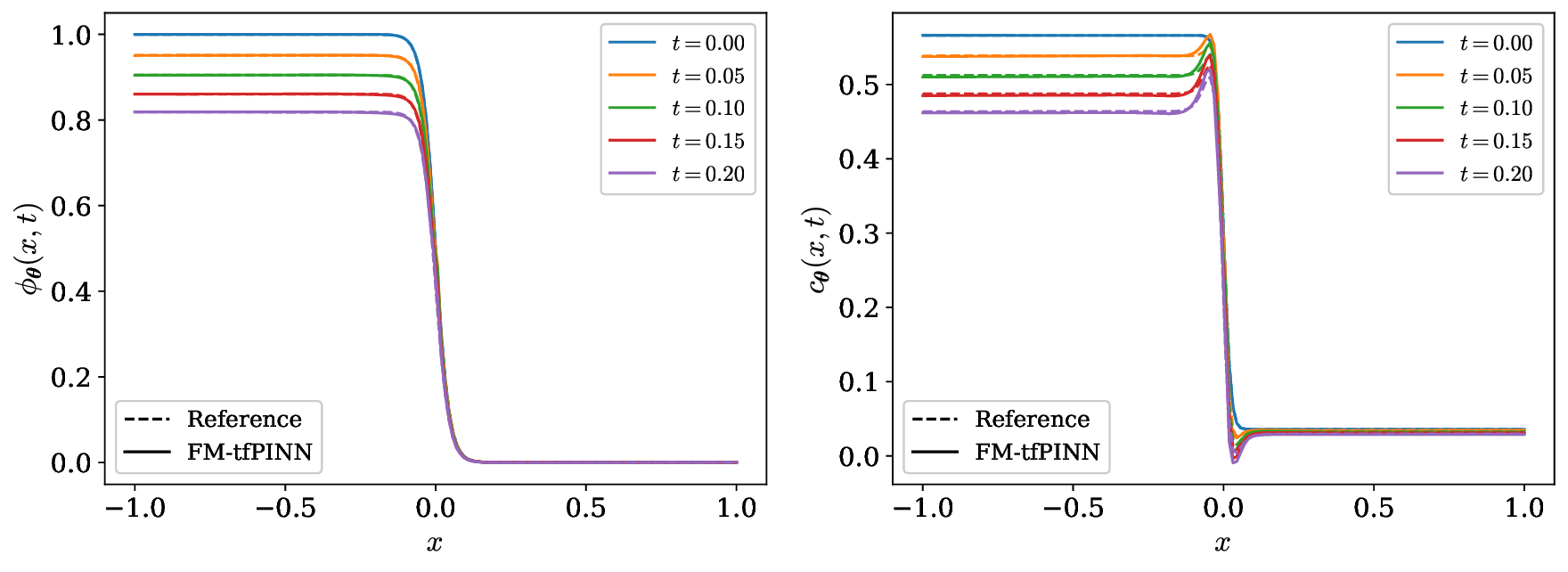}
  \caption{Post-identification field validation for Example~\ref{subsec:ex4_inverse_identification}. The reconstructed phase-field profile $\phi_{\bth}(x,t)$ and concentration profile $c_{\bth}(x,t)$ are compared with the reference solution at representative time levels.}
  \label{fig:ex4_field_profiles}
\end{figure}

Figure~\ref{fig:ex4_field_profiles} provides a post-training validation of the recovered parameters at representative time levels. 
The reconstructed phase field $\phi_{\bth}$ agrees closely with the reference profiles and accurately captures the diffuse corrosion front. 
The concentration field $c_{\bth}$ also follows the reference behavior, including the interfacial redistribution induced by the corrosion process. 

Overall, this inverse experiment shows that FM-tfPINN can recover hidden mobility parameters from limited physical corrosion observations while preserving accurate field reconstruction.
The result extends the forward predictive capability demonstrated in the preceding examples to an inverse setting, where the trainable neural representation and the unknown corrosion parameters are identified simultaneously under the same tempered fractional-memory residual.

\section{Conclusion and Future Perspectives}\label{sect:concl}
In this work, we have developed FM-tfPINN as a fractional-memory generated physics-informed neural framework for tempered time-fractional coupled phase-field systems. The central idea has been to construct the neural approximation through latent memory-source fields whose tempered fractional integrals generate the coupled variables $\bu_{\bth}=(\phi_{\bth},c_{\bth})^\top$.

The framework has combined four complementary components in a unified manner. 
The first component is the tempered fractional-memory neural representation, which embeds the initial state and generates the subsequent evolution through $\mathcal{I}_t^{\alpha,\lambda}$.
The second component is the fast SOE-accelerated shifted memory residual on graded temporal levels $t_{n+\sigma}$, which enables efficient enforcement of the tempered Caputo history contribution. The third component is the interface-aware and residual-adaptive collocation strategy, which concentrates learning effort near the diffuse interface and in regions with larger coupled residuals. The fourth component is the physics-informed objective, which incorporates the coupled residual losses, boundary constraints, admissibility penalties, and, in the inverse setting, sparse physical observation losses. 

The effectiveness of FM-tfPINN has been demonstrated through representative forward and inverse experiments for corrosion-driven phase-field dynamics.
The numerical results show that the proposed framework can resolve coupled phase-field and concentration evolution, capture interface-sensitive corrosion observables, and recover geometry-dependent pitting quantities within a unified memory-generated learning structure. 
The inverse study further indicates the potential of the method for identifying unknown mobility parameters from sparse physical observations, such as $d(t)$, $c_{\Gamma}(t)$, and $c_{\max}(t)$, when full-field measurements are unavailable. 
Overall, these results suggest that FM-tfPINN provides a memory-consistent, interface-aware, and physically interpretable approach for forward prediction and inverse identification in tTFCP systems. 

Future work will focus on extending the proposed framework to more complex coupled phase-field configurations, including three-dimensional pitting, multi-ion electrochemical transport, and experimentally measured corrosion data. Further developments will also consider adaptive loss balancing, uncertainty quantification for noisy and sparse physical observations, and more efficient long-time memory propagation for large-scale tempered fractional phase-field systems.

\section*{Declarations}

\subsection*{Author Contribution}
Shubham Kumar: Methodology, software implementation, numerical experiments,
visualization, and writing -- original draft.

Himanshu Kumar Dwivedi: Conceptualization, software implementation, metho\-dology, mathematical analysis,
validation, investigation, and writing -- original draft.

Matthias Ehrhardt: Supervision, validation, interpretation of results,
writing -- review and editing, and correspondence.

Rajeev: Supervision, validation, project administration, funding acquisition
and writing--review and editing.

All authors have read and approved the final version of the manuscript.
\subsection*{Funding Declarations}
This research has been supported by the Anusandhan National Research Foundation (ANRF), Government of India, under Sanction Letter No.~ANRF/IRG/2025/000423/MS.


\end{document}